\UseRawInputEncoding
\documentclass[12pt]{article}
\usepackage{amssymb}
\usepackage{amsthm}
\usepackage{amsthm,amsmath,indentfirst}
\usepackage{pstricks}
\usepackage{algorithm} 
\usepackage{algorithmic} 
\usepackage{blindtext}
\usepackage{graphicx}
\usepackage{cite}
\usepackage{float}
\usepackage[scriptsize,nooneline]{subfigure}
\usepackage{epstopdf}
\usepackage{txfonts}
\usepackage{amssymb}
\usepackage{amsfonts}
\usepackage[colorlinks,
linkcolor=blue,
anchorcolor=blue,
citecolor=blue
]{hyperref}

\newcommand{\R}{{\mathbb R}}

{}
{}
{}

\begin{document}

\title{Mixed-norm Amalgam Spaces
\thanks{The research was supported by the National Natural Science Foundation of China(12061069).}
}
\author{Houkun Zhang\thanks{ Author E-mail address: zhanghkmath@163.com.},\quad Jiang Zhou\thanks{ Corresponding author E-mail address: zhoujiang@xju.edu.cn.}
\\\\[.5cm]
\small  College of Mathematics and System Sciences, Xinjiang University, Urumqi 830046\\
\small People's Republic of China
}
\date{}
\maketitle
{\bf Abstract:} We introduce the mixed-norm amalgam spaces $(L^{\vec{p}},L^{\vec{s}})(\R^n)$ and $(L^{\vec{p}},L^{\vec{s}})^{\alpha}(\R^n)$, and show their some basic properties. In addition, we find the predual $\mathcal{H}(\vec{p}',\vec{s}\,',\alpha')$ of mixed-norm amalgam spaces $(L^{\vec{p}},\ell^{\vec{s}})^{\alpha}(\R^n)$ by the dual spaces $(L^{\vec{p}'},\ell^{\vec{s}\,'})(\R^n)$ of $(L^{\vec{p}},\ell^{\vec{s}})(\R^n)$, where $(L^{\vec{p}},L^{\vec{s}})(\R^n)=(L^{\vec{p}},\ell^{\vec{s}})(\R^n)$ and $(L^{\vec{p}},L^{\vec{s}})^{\alpha}(\R^n)=(L^{\vec{p}},\ell^{\vec{s}})^{\alpha}(\R^n)$.  Then, we study the strong-type estimates for fractional integral operators $I_{\gamma}$ on mixed-norm amalgam spaces $(L^{\vec{p}},L^{\vec{s}})^{\alpha}(\R^n)$. And, the strong-type estimates of linear commutators $[b,I_{\gamma}]$ generated by $b\in BMO(\R^n)$ and $I_{\gamma}$ on mixed-norm amalgam spaces $(L^{\vec{p}},L^{\vec{s}})^{\alpha}(\R^n)$ are established as well. Furthermore, based on the dual theorem, the characterization of $BMO(\R^n)$ by the boundedness of $[b,I_\gamma]$ from $(L^{\vec{p}},L^{\vec{s}})^{\alpha}(\R^n)$ to $(L^{\vec{q}},L^{\vec{s}})^{\beta}(\R^n)$ is given, which is a new result even for the classical amalgam spaces.
\par
{\bf Keywords:} Mixed norm; Amalgam spaces; Predual; Fractional integral operators; Commutators

{\bf MSC(2000) subject classification:}  42B25; 42B20.

\maketitle

\section{Introduction}\label{sec1}

\par
In fact, mixed-norm Lebesgue spaces, as natural generalizations of the classical Lebesgue spaces $L^p(\R^n)(0<p<\infty)$, were first introduced by Benedek and Panzone \cite{11}. Due to the more precise structure of mixed-norm function spaces than the corresponding classical function spaces, mixed-norm function spaces are of extensive applications in the partial differential equations \cite{7,9,10}. So the mixed-norm function spaces are widely introduced and studied, such as mixed-norm Lorentz spaces \cite{12}, mixed-norm Lorentz-Marcinkiewicz spaces \cite{13}, mixed-norm Orlicz spaces \cite{14}, anisotropic mixed-norm Hardy spaces \cite{15}, mixed-norm Triebel-Lizorkin spaces \cite{16}, mixed Morrey spaces \cite{19,20}, and weak mixed-norm Lebesgue spaces \cite{17}.

The mixed-norm Lebesgue spaces is stated as follows. Let $f$ is a measurable function on $\mathbb{R}^n$ and $0<\vec{p}<\infty$. We say that $f$ belongs to the mixed-norm Lebesgue spaces $L^{\vec{p}}(\mathbb{R}^n)$, if the norm
$$\left\|f\right\|_{L^{\vec{p}}(\mathbb{R}^n)}=\left(\int_{\mathbb{R}}\cdots\left(\int_{\mathbb{R}}\left|f(x)\right|^{p_1}\,dx_1\right)
^{\frac{p_2}{p_1}}\cdots\,dx_n\right)^{\frac{1}{p_n}}<\infty.$$
Note that if $p_1=p_2=\cdots=p_n=p$, then $L^{\vec{p}}(\mathbb{R}^n)$ are reduced to classical Lebesgue spaces $L^p$ and
$$\left\|f\right\|_{L^{\vec{p}}(\mathbb{R}^n)}=\left(\int_{\mathbb{R}^n}\left|f(x)\right|^{p} dx\right)^{\frac{1}{p}}.$$

In this paper, we introduce two new mixed-norm function spaces, mixed-norm amalgam spaces $(L^{\vec{p}},L^{\vec{s}})(\R^n)$ and $(L^{\vec{p}},L^{\vec{s}})^{\alpha}(\R^n)$. Let us to recall some information of classical amalgam spaces.

The amalgam spaces $(L^p,\ell^s)(\R^n)$ were first introduced by Wiener \cite{25} in 1926. However, its systematic study goes back to the works of Holland \cite{5}, which studied the Fourier transform on $\R^n$. Besides that, the spaces have been widely studied \cite{6,18,27,28,30}. It is obvious that Lebesgue space $L^{p}(\R^n)$ coincides with the amalgam space $(L^p,\ell^p)(\R^n)$. It is easy to say that for any $r>0$, the dilation operator $St_{r}^{(p)}:f\mapsto r^{-\frac{n}{p}}f(r^{-1}\cdot)$ is isometric on $L^{p}(\R^n)$. Howeve, amalgam spaces do not have this property. If $p\neq s$, there does not exist $\alpha$ such that $\sup_{r>0}\|St_{r}^{(\alpha)}(f)\|_{(L^p,\ell^s)}<\infty$, although $St_{r}^{(\alpha)}(f)\in(L^{p},\ell^s)(\R^n)$ for all $f\in(L^p,\ell^s)(\R^n)$, $\rho>0$ and $\alpha>0$ \cite{8}. The amalgam spaces $(L^p,\ell^s)^\alpha(\R^n)$ compensate this shortcomings. The functions spaces $(L^p,\ell^s)^\alpha(\R^n)$ were introduced by Fofana in 1988, which consist of $f\in(L^p,\ell^s)(\R^n)$ and satisfying $\sup_{r>0}\|St_{r}^{(\alpha)}(f)\|_{(L^p,\ell^s)}<\infty$.

The fractional power of the Laplacian operators $\triangle$ are defined by
$$\left((-\triangle)^{\gamma/2}(f)\right)^{\wedge}(\xi)=(2\pi|\xi|)^{\gamma}\hat{f}(\xi).$$
Comparing this to the Fourier transform of $|x|^{-\gamma}$, $0<\gamma<n$, we are led to define the so-called fractional integral operators $I_{\gamma}$ by
$$I_\gamma f(x)=(-\triangle)^{\gamma/2}(f)(x)=C_\gamma\int_{\mathbb{R}^n}\frac{f(y)}{|x-y|^{n-\gamma}}dy,$$
where
$$C_\gamma^{-1}=\frac{\pi^{n/2}2^{\gamma}\Gamma(\gamma/2)}{\Gamma((n-\gamma)/2)}.$$
The fractional integral operators play an important role in harmonic analysis. An important application, via the well-knowing Hardy-Littlewood-Sobolev theorem, is in proving the Sobolev embedding theorem. In this paper, we investigate the generalization of Hardy-Littlewood-Sobolev theorem on mixed-norm amalgam spaces.

The boundedness properties of $I_{\gamma}$ between various function spaces have been studied extensively. In 1960, Benedek and Panzone first study the boundedness of $I_{\gamma}$ from mixed-norm Lebesgue spaces $L^{\vec{p}}(\R^n)$ to mixed-norm Lebesgue spaces $L^{\vec{q}}(\R^n)$ \cite{11}, which is generalization of the classical Hardy-Littlewood-Sobolev theorem (see \cite{1}). In 2021, Zhang and Zhou improve the theorem on mixed-norm Lebesgue spaces, which is stated as follows.

\textbf{Lemma 1.1.} (see \cite{2}) Let $0<\gamma<n$ and $1<\vec{p},\vec{q}<\infty$. Then
$$1<\vec{p}\le\vec{q}<\infty,~\gamma=\sum_{i=1}^n\frac{1}{p_i}-\sum_{i=1}^n\frac{1}{q_i}$$
if and only if
$$\|I_{\gamma}f\|_{L_{\vec{q}}}\lesssim\|f\|_{L_{\vec{p}}}.$$

For a locally integrable function $b$, the commutators of fractional integral operators $I_\gamma$ are defined by
$$[b,I_\gamma]f(x):=b(x)I_\gamma f(x)-I_\gamma(bf)(x)=C_\gamma\int_{\mathbb{R}^n}\frac{(b(x)-b(y))f(y)}{|x-y|^{n-\gamma}}dy,$$
which were introduced by Chanillo in \cite{3}. These commutators also can be used to study theory of Hardy spaces $H^p(\R^n)$\cite{26}. In 2019, Nogayama given an characterization of $BMO(\mathbb{R}^n)$ spaces via the $(\mathcal{M}_{\vec{p}}^{p_0},\mathcal{M}_{\vec{q}}^{q_0})$-boundedness of $[b,I_\gamma]$\cite{20}. In the 2021, the result is improved on mixed-norm Lebesgue in \cite{4}, which is stated as follows.

\textbf{Lemma 1.2.} (see \cite{4}) Let $0<\gamma<n,~1<\vec{p}\le\vec{q}<\infty$ and
$$\gamma=\sum_{i=1}^n\frac{1}{p_i}-\sum_{i=1}^n\frac{1}{q_i}.$$
Then the following conditions are equivalent:\\
(\romannumeral1) $b\in BMO(\mathbb{R}^n)$.\\
(\romannumeral2) $[b,I_\gamma]$ is bounded from $L^{\vec{p}}(\mathbb{R}^n)$ to $L^{\vec{q}}(\mathbb{R}^n)$.

Now, let us recall the definition of $BMO(\R^n)$. $BMO(\R^n)$ is the Banach function space modulo constants with the norm $\|\cdot\|_{BMO}$ defined by
$$\|b\|_{BMO}=\sup_{B\subset\mathbb{R}^n}\frac{1}{|B|}\int_{B}|b(y)-b_B|dy<\infty,$$
where the supremum is taken over all balls $B$ in $\R^n$ and $b_B$ stands for the mean value of $b$ over $B$; that is, $b_B:=(1/|B|)\int_Bb(y)dy$. By John-Nirenberg inequality,
$$\|b\|_{BMO}\sim\sup_{B\subset\mathbb{R}^n}\frac{\|b-b_B\|_{L^p}}{\|\chi_{B}\|_{L^p}},~~1<p<\infty.$$
It is also right if we replace $L^p$-norm by mixed-norm $L^{\vec{p}}$-norm (see Lemma 4.1).

We firstly define mixed-norm amalgam spaces, which can be considered as an extension of classical amalgam spaces. It is natural and important to study the boundedness of $I_\gamma$ and $[b,I_\gamma]$ in these new spaces. Before that, we also study some properties of the new spaces. This paper is organized as follows. In Section 2, we state definitions of mixed-norm amalgam spaces, some properties of mixed-norm amalgam spaces, and the main results of the present paper. We will give the proof of some properties of mixed-norm amalgam spaces in Section 3. The predual of mixed-norm amalgam spaces is studied in Section 4. In the Section 5 and 6, we prove the boundedness of $I_\gamma$ and their commutators generated by $b\in BMO(\R^n)$. In the final section, We study the necessary condition of the boundedness of $[b,I_\gamma]$ from $(L^{\vec{p}},L^{\vec{s}})^\alpha(\R^n)$ to $(L^{\vec{q}},L^{\vec{s}})^\beta(\R^n)$, which is a new result even for the classical amalgam spaces.

Next, we make some conventions and recall some notions. Let $\vec{p}=(p_1,p_2,$ $\cdots,p_n),~\vec{q}=(q_1,q_2,\cdots,q_n)$, $\vec{s}=(s_1,s_2,\cdots,s_n)$, are n-tuples and $1<p_i,q_i,s_i<\infty,~i=1,2,\cdots,n$. We define that  if $\varphi(a,b)$ is a relation or equation among numbers, $\varphi(\vec{p},\vec{q})$ will mean that $\varphi(p_i,q_i)$ hords for each $i$. For example, $\vec{p}<\vec{q}$ means that $p_i<q_i$ holds for each $i$ and $\frac{1}{\vec{p}}+\frac{1}{\vec{p}\,'}=1$ means $\frac{1}{p_i}+\frac{1}{p_i'}=1$ hold for each $i$. The symbol $B$ denote the open ball and $B(x,r)$ denote the open ball centered at $x$ of radius $r$. Let $\rho B(x,r)=B(x,\rho r)$, where $\rho>0$. Let $L^{\vec{p}}=L^{\vec{p}}(\R^n)$, $(L^{\vec{p}},L^{\vec{s}})=(L^{\vec{p}},L^{\vec{s}})(\R^n)$ and $(L^{\vec{p}},L^{\vec{s}})^\alpha=(L^{\vec{p}},L^{\vec{s}})^\alpha(\R^n)$. $A\sim B$ means that $A$ is equivalent to $B$, that is, $A\lesssim B(A\le CB)$ and $B\lesssim A(B\le CA)$, where $C$ is a positive constant. Through all paper, each positive constant $C$ is not necessarily equal.

\section{Mixed-norm amalgam spaces $(L^{\vec{p}},L^{\vec{s}})(\R^n)$ and $(L^{\vec{p}},L^{\vec{s}})^\alpha(\R^n)$}\label{sec2}
\par
In this section, we firstly present the definitions of mixed-norm amalgam spaces and some properties of mixed-norm amalgam spaces in Section 2.1, and then main theorems are showed in Section 2.2.
\subsection{Definitions and properties}
In this section, we present the definitions of mixed-norm amalgam spaces $(L^{\vec{p}},L^{\vec{s}})(\R^n)$ and $(L^{\vec{p}},L^{\vec{s}})^\alpha(\R^n)$ and their properties. Firstly, the definitions of mixed-norm amalgam spaces $(L^{\vec{p}},L^{\vec{s}})(\R^n)$ and $(L^{\vec{p}},L^{\vec{s}})^\alpha(\R^n)$ are given as follows.

\textbf{Definition 2.1} Let $1\le\vec{p},\vec{s},\alpha\le\infty$. We define two types of amalgam spaces of $L^{\vec{p}}(\mathbb{R}^n)$ and $L^{\vec{s}}(\mathbb{R}^n)$. If measurable functions $f$ satisfy $f\in L^{1}_{loc}(\mathbb{R}^n)$, then
$$(L^{\vec{p}},L^{\vec{s}})(\mathbb{R}^n) :=\left\{f:\|f\|_{(L^{\vec{p}},L^{\vec{s}})}<\infty\right\}$$
and
$$(L^{\vec{p}},L^{\vec{s}})^{\alpha}(\mathbb{R}^n) :=\left\{f:\|f\|_{(L^{\vec{p}},L^{\vec{s}})^{\alpha}}<\infty\right\},$$
where
$$\|f\|_{(L^{\vec{p}},L^{\vec{s}})}=\left\|\|f\chi_{B(\cdot,1)}\|_{L^{\vec{p}}}\right\|_{L^{\vec{s}}} =\left(\int_{\mathbb{R}}\cdots\left(\int_{\mathbb{R}}\|f\chi_{B(y,1)}\|_{L^{\vec{p}}}^{s_1}\,dy_1\right)
^{\frac{s_2}{s_1}}\cdots\,dy_n\right)^{\frac{1}{s_n}}$$
and
\begin{align*}
\|f\|_{(L^{\vec{p}},L^{\vec{s}})^\alpha}
&=\sup_{r>0}\left\||B(\cdot,r)|^{\frac{1}{\alpha}-\frac{1}{n}\sum_{i=1}^{n}\frac{1}{p_i}-\frac{1}{n}\sum_{i=1}^{n}\frac{1}{s_i}} \|f\chi_{B(\cdot,r)}\|_{L^{\vec{p}}(\mathbb{R}^n)}\right\|_{L^{\vec{s}}(\mathbb{R}^n)}\\
&=\sup_{r>0}\left(\int_{\mathbb{R}}\cdots\left(\int_{\mathbb{R}} \left(|B(y,r)|^{\frac{1}{\alpha}-\frac{1}{n}\sum_{i=1}^{n}\frac{1}{p_i}-\frac{1}{n}\sum_{i=1}^{n}\frac{1}{s_i}} \|f\chi_{B(y,r)}\|_{L^{\vec{p}}(\mathbb{R}^n)}\right)^{s_1}\,dy_1\right)
^{\frac{s_2}{s_1}}\cdots\,dy_n\right)^{\frac{1}{s_n}}
\end{align*}
with the usual modification for $p_i=\infty$ or $s_i=\infty$.

Next, we claim that the mixed-norm amalgam spaces defined in Definition 2.1 are Banach spaces.

\textbf{Proposition 2.2.} Let $1\le\vec{p},\vec{s},\alpha\le\infty$. Mixed norm amalgam spaces $(L^{\vec{p}},L^{\vec{s}})(\R^n)$ and $(L^{\vec{p}},L^{\vec{s}})^\alpha(\R^n)$ are also Banach spaces.

The following proposition shows that the necessary relationship of the index $\vec{p},\vec{s}$ and $\alpha$.

\textbf{Proposition 2.3.} The spaces $(L^{\vec{p}},L^{\vec{s}})^{\alpha}(\R^n)$ are nontrivial if and only if $\frac{1}{n}\sum_{i=1}^n\frac{1}{s_i}\leq\frac{1}{\alpha}\leq\frac{1}{n}\sum^{n}_{i=1}\frac{1}{p_{i}}$.

By Definition 2.1, if $p_i=p$ and $s_i=s$ for each $i$, then
$$(L^{\vec{p}},L^{\vec{s}})(\mathbb{R}^n)=(L^p,L^s)(\mathbb{R}^n), ~(L^{\vec{p}},L^{\vec{s}})^{\alpha}(\mathbb{R}^n)=(L^p,L^s)^{\alpha}(\mathbb{R}^n).$$
In particular, If $s_i=\infty$ for each $i$ and $\frac{1}{\alpha}\leq\frac{1}{n}\sum^{n}_{i=1}\frac{1}{p_{i}}$, then
$$(L^{\vec{p}},L^{\vec{s}})^{\alpha}(\mathbb{R}^n)=\mathcal{M}_{\vec{p}}^{\alpha}(\mathbb{R}^n),$$
where $\mathcal{M}_{\vec{p}}^{\alpha}(\mathbb{R}^n)$ is mixed Morrey spaces defined as \cite{19,20}. Finally, we study the relationship between the mixed-norm amalgam spaces.

\textbf{Proposition 2.4.} Let $1\le\vec{p},\vec{q},\vec{s}\le\infty$, $\frac{1}{n}\sum_{i=1}^n\frac{1}{s_i}\leq\frac{1}{\alpha}\leq\frac{1}{n}\sum^{n}_{i=1}\frac{1}{p_{i}}$, and $\frac{1}{n}\sum_{i=1}^n\frac{1}{s_i}\leq\frac{1}{\alpha}\leq\frac{1}{n}\sum^{n}_{i=1}\frac{1}{q_{i}}$. Then\\
(\romannumeral1) $(L^{\vec{p}},L^{\vec{s}})^\alpha(\R^n)\subset(L^{\vec{p}},L^{\vec{s}})(\R^n)$ with $\|f\|_{(L^{\vec{p}},L^{\vec{s}})}\le\|f\|_{(L^{\vec{p}},L^{\vec{s}})^{\alpha}}$;\\
(\romannumeral2) If $\vec{p}\le\vec{q}$, $(L^{\vec{q}},L^{\vec{s}})^{\alpha}(\R^n)\subseteq(L^{\vec{p}},L^{\vec{s}})^{\alpha}(\R^n)$ with $\|f\|_{(L^{\vec{p}},L^{\vec{s}})^{\alpha}}\le\|f\|_{(L^{\vec{q}},L^{\vec{s}})^{\alpha}}$.

\subsection{Main Theorems}
\par
In this section, we show the main theorems in this paper. Before all, we give the equivalent norms of mixed-norm amalgam spaces. Let $Q_{r,k}=r[k+[0,1)^n]$ and
$$\left\|\{a_k\}_{k\in \mathbb{Z}^n}\right\|_{\ell^{\vec{s}}} :=\left(\sum_{k_n\in\mathbb{Z}}\cdots
\left(\sum_{k_1\in\mathbb{Z}}|a_k|^{s_1}\right)^{\frac{s_{2}}{s_1}}\cdots\right)^{\frac{1}{s_n}}$$
with the usual modification for $s_i=\infty$.

\textbf{Proposition 2.5.} Let $1\le\vec{p},\vec{s},\alpha\le\infty$ and $\frac{1}{n}\sum_{i=1}^{n}\frac{1}{s_i}\leq\frac{1}{\alpha}\leq\frac{1}{n}\sum^{n}_{i=1}\frac{1}{p_{i}}$. We define two types ``discrete" mixed-norm amalgam spaces.
$$(L^{\vec{p}},\ell^{\vec{s}})(\R^n):=\left\{f\in L_{loc}^{1}:\|f\|_{\vec{p},\vec{s}} :=\left\|\left\{\|f\chi_{Q_{1,k}}\|_{\vec{p}}\right\}_{k\in\mathbb{Z}^n}\right\|_{\ell^{\vec{s}}}<\infty\right\}$$
and
$$(L^{\vec{p}},\ell^{\vec{s}})^\alpha(\R^n):=\left\{f\in L_{loc}^{1}:\|f\|_{\vec{p},\vec{s},\alpha} :=\sup_{r>0}r^{\frac{n}{\alpha}-\sum_{i=1}^n\frac{1}{p_i}}{_r\|f\|_{\vec{p},\vec{s}}}<\infty\right\},$$
where
$$_r\|f\|_{\vec{p},\vec{s}}:=\left\|\left\{\|f\chi_{Q_{r,k}}\|_{\vec{p}}\right\}_{k\in\mathbb{Z}^n}\right\|_{\ell^{\vec{s}}}.$$
In fact, we have
$$(L^{\vec{p}},L^{\vec{s}})(\mathbb{R}^n)=(L^{\vec{p}},\ell^{\vec{s}})(\mathbb{R}^n)\text{ and }(L^{\vec{p}},L^{\vec{s}})^\alpha(\mathbb{R}^n)=(L^{\vec{p}},\ell^{\vec{s}})^\alpha(\mathbb{R}^n).$$

According to Proposition 2.5, we give the definition of the predual of mixed-norm amalgam spaces $(L^{\vec{p}},\ell^{\vec{s}})^\alpha(\mathbb{R}^n)$.

\textbf{Definition 2.6.} Let $1\le\vec{p},\vec{s},\alpha\le\infty$ and $\frac{1}{n}\sum_{i=1}^{n}\frac{1}{s_i}\leq\frac{1}{\alpha}\leq\frac{1}{n}\sum^{n}_{i=1}\frac{1}{p_{i}}$. The space $\mathcal{H}(\vec{p}',\vec{s}\,',\alpha')$ is defined as the set of all elements of $L^1_{loc}(\R^n)$ for which there exist a sequence $\{(c_j,r_j,f_j)\}_{j\ge 1}$ of elements of $\mathbb{C}\times(0,\infty)\times(L^{\vec{p}'},\ell^{\vec{s}\,'})(\R^n)$ such that\\
$$f:=\sum_{j\ge 1}c_j St_{r_j}^{(\alpha')}(f_j)\text{ in the sense of }L^1_{loc}(\R^n);\eqno{(2.1)}$$
$$\|f_j\|_{\vec{p}',\vec{s}\,'}\le 1,j\ge 1; \eqno{(2.2)}$$
$$\sum_{j\le 1}|c_j|<\infty.\eqno{(2.3)}$$
We will always refer to any sequence $\{(c_j,r_j,f_j)\}_{j\ge 1}$ of elements of $\mathbb{C}\times(0,\infty)\times(L^{\vec{p}'},\ell^{\vec{s}\,'})(\R^n)$ satisfying (2.1)-(2.3) as block decomposition of $f$. For any element $f$ of $\mathcal{H}(\vec{p}',\vec{s}\,',\alpha')$, we set
$$\|f\|_{\mathcal{H}(\vec{p}',\vec{s}\,',\alpha')}:=\inf\left\{\sum_{j\ge 1}|c_j|:f:=\sum_{j\ge 1}c_j St_{r_j}^{(\alpha')}f_j\right\},$$
where the infimum is taken over all block decomposition of $f$.

\textbf{Theorem 2.7.} (\romannumeral1) Let $1\le\vec{p},\vec{s}\le\infty$, and $\frac{1}{n}\sum_{i=1}^n\frac{1}{s_i}\leq\frac{1}{\alpha}\leq\frac{1}{n}\sum^{n}_{i=1}\frac{1}{p_{i}}$. If $g\in(L^{\vec{p}},\ell^{\vec{s}})^\alpha$ and $f\in\mathcal{H}(\vec{p}',\vec{s}\,',\alpha')$, we obtain $fg\in L^1(\R^n)$ and
$$\left|\int_{\R^n}f(x)g(x)dx\right|\le\|g\|_{(L^{\vec{p}},\ell^{\vec{s}})^\alpha} \|f\|_{\mathcal{H}(\vec{p}',\vec{s}\,',\alpha')}.\eqno{(2.4)}$$

(\romannumeral2) Let $1<\vec{p},\vec{s}\le\infty$ and $\frac{1}{n}\sum_{i=1}^n\frac{1}{s_i}\leq\frac{1}{\alpha}\leq\frac{1}{n}\sum^{n}_{i=1}\frac{1}{p_{i}}$. The operator $T:g\mapsto T_g$ defined as
$$<T_g,f>=\int_{\R^n}f(x)g(x)dx,~~g\in(L^{\vec{p}},\ell^{\vec{s}})^{\alpha}(\R^n)\text{ and }f\in\mathcal{H}{(\vec{p}',\vec{s}\,',\alpha')}$$
is an isomestric isomorphism of $(L^{\vec{p}},\ell^{\vec{s}})^{\alpha}(\R^n)$ into $\mathcal{H}{(\vec{p}',\vec{s}\,',\alpha')}^{*}$.

Now, We show the boundedness of fractional integral operators on mixed-norm amalgam spaces.

\textbf{Theorem 2.8.} Let $0<\gamma<n$, $1<\vec{p},\vec{q}<\infty$, $1<\vec{s}\le\infty$,
 $\frac{1}{n}\sum_{i=1}^{n}\frac{1}{s_i}\le\frac{1}{\alpha}\le\frac{1}{n}\sum_{i=1}^n\frac{1}{p_i}$, and $\frac{1}{n}\sum_{i=1}^{n}\frac{1}{s_i}\le\frac{1}{\beta}\le\frac{1}{n}\sum_{i=1}^n\frac{1}{q_i}$. Assume that  $\gamma=\sum_{i=1}^n\frac{1}{p_i}-\sum_{i=1}^{n}\frac{1}{q_i}$. Then the fractional integral operators $I_\gamma$ are bounded from $(L^{\vec{p}},L^{\vec{s}})^{\alpha}(\R^n)$ to $(L^{\vec{q}},L^{\vec{s}})^{\beta}(\R^n)$ if and only if $$\gamma=\frac{n}{\alpha}-\frac{n}{\beta}.$$

\textbf{Remark 2.9.} In fact, the condition $\gamma=\frac{n}{\alpha}-\frac{n}{\beta}$ is necessary for the boundedness of fractional integral operators $I_\gamma$. Let $\delta_tf(x)=f(tx)$, where $(t>0)$. Then,
$$I_{\gamma}(\delta_t f)=t^{-\gamma}\delta_t I_{\gamma}(f)$$
$$\|\delta_{t^{-1}}f\|_{_{(L^{\vec{q}},L^{\vec{s}})^{\beta}}} =t^{\frac{n}{\beta}-\frac{1}{n}\sum_{i=1}^n\frac{1}{s_i}}\|f\|_{_{(L^{\vec{q}},L^{\vec{s}})^{\beta}}}.$$
$$\|\delta_{t}f\|_{_{(L^{\vec{p}},L^{\vec{r}})^{\alpha}}} =t^{-\frac{n}{\alpha}+\sum_{i=1}^{n}\frac{1}{r_i}}\|f\|_{_{(L^{\vec{p}},L^{\vec{r}})^{\alpha}}}.$$
Thus, by the boundedness of $I_\gamma$ from $(L^{\vec{p}},L^{\vec{r}})^{\alpha}(\R^n)$ to $(L^{\vec{q}},L^{\vec{s}})^{\beta}(\R^n)$,
\begin{align*}
\|I_{\gamma}f\|_{(L^{\vec{q}},L^{\vec{s}})^{\beta}}
&=t^{\gamma}\|\delta_{t^{-1}}I_{\gamma}(\delta_tf)\|_{(L^{\vec{q}},L^{\vec{s}})^{\beta}}\\
&=t^{\gamma+\frac{n}{\beta}-\sum_{i=1}^n\frac{1}{s_i}}\|I_{\gamma}(\delta_tf)\|_{(L^{\vec{q}},L^{\vec{s}})^{\beta}}\\
&\lesssim t^{\gamma+\frac{n}{\beta}-\sum_{i=1}^n\frac{1}{s_i}}\|\delta_tf\|_{(L^{\vec{p}},L^{\vec{r}})^{\alpha}}\\
&=t^{\gamma+\frac{n}{\beta}-\sum_{i=1}^n\frac{1}{s_i}-\frac{n}{\alpha}+\sum_{i=1}^n\frac{1}{r_i}}
\|f\|_{(L^{\vec{p}},L^{\vec{r}})^{\alpha}}.
\end{align*}
Thus, $\gamma=\frac{n}{\alpha}-\sum_{i=1}^{n}\frac{1}{r_i}-\frac{n}{\beta}+\sum_{i=1}^{n}\frac{1}{s_i}$ and $\gamma=\frac{n}{\alpha}-\frac{n}{\beta}$ when $\vec{s}=\vec{r}$.

Let $[b,I_\gamma]$ be the linear commutators generated by $I_\gamma$ and $BMO$ function $b$. For the strong-type estimates of $[b,I_\gamma]$ on the mixed-norm amalgam spaces, we have the following result.

\textbf{Theorem 2.10.} Let $0<\gamma<n$, $1<\vec{p},\vec{q}<\infty$, $1<\vec{s}\le\infty$, $\frac{1}{n}\sum_{i=1}^n\frac{1}{s_i}\le\frac{1}{\alpha}\le\frac{1}{n}\sum_{i=1}^n\frac{1}{p_i}$, and $\frac{1}{n}\sum_{i=1}^n\frac{1}{s_i}\le\frac{1}{\beta}\le\frac{1}{n}\sum_{i=1}^n\frac{1}{q_i}$. Assume that $\gamma=\sum_{i=1}^n\frac{1}{p_i}-\sum_{i=1}^{n}\frac{1}{q_i}=\frac{n}{\alpha}-\frac{n}{\beta}$. If $b\in BMO(\R^n)$, then the linear commutators $[b,I_\gamma]$ are bounded from $(L^{\vec{p}},L^{\vec{s}})^{\alpha}(\R^n)$ to $(L^{\vec{q}},L^{\vec{s}})^{\beta}(\R^n)$.

In fact, if the linear commutators $[b,I_\gamma]$ are bounded from $(L^{\vec{p}},L^{\vec{s}})^{\alpha}(\R^n)$ to $(L^{\vec{q}},L^{\vec{s}})^{\beta}(\R^n)$, then $b\in BMO(\R^n)$. This result can be stated as follows.

\textbf{Theorem 2.11.}  Let $0<\gamma<n$, $1<\vec{p},\vec{q}<\infty$, $1<\vec{s}\le\infty$, $\frac{1}{n}\sum_{i=1}^n\frac{1}{s_i}\le\frac{1}{\alpha}\le\frac{1}{n}\sum_{i=1}^n\frac{1}{p_i}$, and $\frac{1}{n}\sum_{i=1}^n\frac{1}{s_i}\le\frac{1}{\beta}\le\frac{1}{n}\sum_{i=1}^n\frac{1}{q_i}$. Assume that $\gamma=\sum_{i=1}^n\frac{1}{p_i}-\sum_{i=1}^{n}\frac{1}{q_i}=\frac{n}{\alpha}-\frac{n}{\beta}$. If the linear commutators $[b,I_\gamma]$ are bounded from $(L^{\vec{p}},L^{\vec{s}})^{\alpha}(\R^n)$ to $(L^{\vec{q}},L^{\vec{s}})^{\beta}(\R^n)$, then $b\in BMO(\R^n)$.

Theorem 2.11 is proved by Proposition 2.5 and Theorem 2.7. By this new result, we can get the following result.

\textbf{Corollary 2.12.}  Let $0<\gamma<n$, $1<\vec{p},\vec{q}<\infty$, $1<\vec{s}\le\infty$, $\frac{1}{n}\sum_{i=1}^n\frac{1}{s_i}\le\frac{1}{\alpha}\le\frac{1}{n}\sum_{i=1}^n\frac{1}{p_i}$, and $\frac{1}{n}\sum_{i=1}^n\frac{1}{s_i}\le\frac{1}{\beta}\le\frac{1}{n}\sum_{i=1}^n\frac{1}{q_i}$. If $\gamma=\sum_{i=1}^n\frac{1}{p_i}-\sum_{i=1}^{n}\frac{1}{q_i}=\frac{n}{\alpha}-\frac{n}{\beta}$, then the following statements are equivalent:\\
(\romannumeral1) The linear commutators $[b,I_\gamma]$ are bounded from $(L^{\vec{p}},L^{\vec{s}})^{\alpha}(\R^n)$ to $(L^{\vec{q}},L^{\vec{s}})^{\beta}(\R^n)$;\\
(\romannumeral2) $b\in BMO(\R^n)$.

\section{The Proofs of Proposition 2.2-2.5}
\par
In this section, we give the proofs of properties of mixed-norm amalgam spaces.

\textbf{Proof of Proposition 2.2.} First, we will check the triangle inequality. For $f,g\in(L^{\vec{p}},L^{\vec{s}})^{\alpha}(\mathbb{R}^n)$,
\begin{align*}
\|f+g\|_{(L^{\vec{p}},L^{\vec{s}})^\alpha}
&=\sup_{r>0}\left\||B(\cdot,r)|^{\frac{1}{\alpha}-\frac{1}{n}\sum_{i=1}^{n}\frac{1}{p_i}-\frac{1}{n}\sum_{i=1}^n\frac{1}{s_i}} \|(f+g)\chi_{B(\cdot,r)}\|_{L^{\vec{p}}}\right\|_{L^{\vec{s}}}\\
&\le\sup_{r>0}\left\||B(\cdot,r)|^{\frac{1}{\alpha}-\frac{1}{n}\sum_{i=1}^{n}\frac{1}{p_i}-\frac{1}{n}\sum_{i=1}^n\frac{1}{s_i}} \|f\chi_{B(\cdot,r)}\|_{L^{\vec{p}}}\right\|_{L^{\vec{s}}}\\ &+\sup_{r>0}\left\||B(\cdot,r)|^{\frac{1}{\alpha}-\frac{1}{n}\sum_{i=1}^{n}\frac{1}{p_i}-\frac{1}{n}\sum_{i=1}^n\frac{1}{s_i}} \|g\chi_{B(\cdot,r)}\|_{L^{\vec{p}}}\right\|_{L^{\vec{s}}}\\
&=\|f\|_{(L^{\vec{p}},L^{\vec{s}})^\alpha}+\|g\|_{(L^{\vec{p}},L^{\vec{s}})^\alpha}.
\end{align*}
The positivity and the homogeneity are both clear. Thus, we prove that $(L^{\vec{p}},L^{\vec{s}})^\alpha(\mathbb{R}^n)$ are spaces with norm $\|\cdot\|_{(L^{\vec{p}},L^{\vec{s}})^\alpha}$. It remains to check the completeness. Without loss the generality, let a Cauchy sequence $\{f_j\}_{j=1}^{\infty}\subset (L^{\vec{p}},L^{\vec{s}})^\alpha(\mathbb{R}^n)$ satisfy $$\|f_{j+1}-f_j\|_{(L^{\vec{p}},L^{\vec{s}})^\alpha}<2^{-j}.$$
We write $f=f_1+\sum_{j=1}^{\infty}(f_{j+1}-f_j)=\lim_{j\rightarrow\infty}f_j$. Then,
$$\|f\|_{(L^{\vec{p}},L^{\vec{s}})^{\alpha}}\le \|f_{1}\|_{(L^{\vec{p}},L^{\vec{s}})^{\alpha}}
+\sum^{\infty}_{j=1}\|f_{j+1}-f_{j}\|_{(L^{\vec{p}},L^{\vec{s}})^{\alpha}}<\infty.$$
Thus, for almost everywhere $x\in\mathbb{R}^n$,
$$f(x)=f_1(x)+\sum_{j=1}^{\infty}(f_{j+1}(x)-f_j(x))\le |f_1(x)|+\sum_{j=1}^{\infty}|f_{j+1}(x)-f_j(x)|<\infty$$
and $f\in(L^{\vec{p}},L^{\vec{s}})^\alpha(\mathbb{R}^n)$. Furthermore,
\begin{align*}
\|f-f_{J}\|_{(L^{\vec{p}},L^{\vec{s}})^{\alpha}}
&=\|\sum^{\infty}_{j=1}(f_{j+1}-f_{j})-\sum^{J-1}_{j=1}(f_{j+1}-f_{j})\|_{(L^{\vec{p}},L^{\vec{s}})^{\alpha}}\\
&\le\sum^{\infty}_{j=J}\|f_{j+1}-f_{j}\|_{(L^{\vec{p}},L^{\vec{s}})^{\alpha}}\\
&\le 2\cdot 2^{-J}
\end{align*}
and
$$\lim_{J\rightarrow\infty}\|f-f_{J}\|_{(L^{\vec{p}},L^{\vec{s}})^{\alpha}}=0.$$
So, we prove that $(L^{\vec{p}},L^{\vec{s}})^{\alpha}$ are Banach spaces. By the same discussion, we can prove $(L^{\vec{p}},L^{\vec{s}})(\R^n)$ also are Banach spaces. $~~~~\blacksquare$

\textbf{Proof of Proposition 2.3.} We prove these by contradiction. In fact, by Lebesgue differential theorem in the mixed-norm Lebesgue spaces \cite{2}, we know
$$\lim_{r\rightarrow 0}\frac{\|f\chi_{B(x,r)}\|_{L^{\vec{p}}}}{\|\chi_{B(x,r)}\|_{L^{\vec{p}}}}=f(x)\text{~a.e.~}x\in\R^n.$$
Thus, if $\frac{1}{n}\sum_{i=1}^n\frac{1}{s_i}-\frac{1}{\alpha}>0$ and $f\neq 0$,
$$ \lim_{r\rightarrow 0}|B(y,r)|^{\frac{1}{\alpha}-\frac{1}{n}\sum_{i=1}^n\frac{1}{s_i}}
\frac{\|f\chi_{B(x,r)}\|_{L^{\vec{p}}}}{\|\chi_{B(x,r)}\|_{L^{\vec{p}}}}=\infty.$$
By this, we prove $\frac{1}{n}\sum_{i=1}^n\frac{1}{s_i}<\frac{1}{\alpha}$.

If $\frac{1}{\alpha}-\frac{1}{n}\sum^{n}_{i=1}\frac{1}{p_{i}}>0$, then we claim $\|\chi_{B(x_{0},r_{0})}\|_{(L^{\vec{p}},L^{\vec{s}})^{\alpha}}=\infty$ for any ball $B(x_0,r_0)$, which show that $(L^{\vec{p}},L^{\vec{s}})^{\alpha}(\R^n)$ are trivial spaces. Hence, we acquire $\frac{1}{\alpha}\leq\frac{1}{n}\sum^{n}_{i=1}\frac{1}{p_{i}}$. Indeed, if $x\in B(x_0,\frac{r}{2})$ and $2r_0<r$, then for any $y\in B(x_0,r_0)$, we have
$$|x-y|\le|x_0-x|+|x_0-y|\le\frac{r}{2}+r_0<r,$$
that is $B(x_0,r_0)\subset B(x,r)$. Therefore,
\begin{align*}
\|\chi_{B(x_{0},r_{0})}\|_{(L^{\vec{p}},L^{\vec{s}})^{\alpha}}
&\sim\sup_{r>0}r^{\frac{n}{\alpha}-\frac{1}{n}\sum_{i=1}^n\frac{1}{s_i}-\sum^{n}_{i=1}\frac{1}{p_{i}}}
\left\| \|\chi_{B(x_{0},r_0)}\chi_{B(\cdot,r)}\|_{L^{\vec{p}}} \right\|_{L^{\vec{s}}}\\
&\ge\sup_{r>2r_0}r^{\frac{n}{\alpha}-\frac{1}{n}\sum_{i=1}^n\frac{1}{s_i}-\sum^{n}_{i=1}\frac{1}{p_{i}}} \left\|\chi_{B(x_0,\frac{r}{2})}\|\chi_{B(x_{0},r_0)}\|_{L^{\vec{p}}} \right\|_{L^{\vec{s}}}\\
&\gtrsim\sup_{r>2r_0} r^{\frac{n}{\alpha}-\frac{1}{n}\sum_{i=1}^n\frac{1}{s_i}-\sum^{n}_{i=1}\frac{1}{p_{i}}}\cdot r^{\frac{1}{n}\sum_{i=1}^n\frac{1}{s_i}}\\
&\geq\lim_{r\rightarrow+\infty} r^{\frac{n}{\alpha}-\sum^{n}_{i=1}\frac{1}{p_{i}}}\\
&=+\infty.
\end{align*}

For the opposite side, it is easy to prove $\chi_{B(0,1)}\in(L^{\vec{p}},L^{\vec{s}})^{\alpha}(\R^n)$ if $\frac{1}{n}\sum_{i=1}^n\frac{1}{s_i}\le\frac{1}{\alpha}\le\frac{1}{n}\sum_{i=1}^n\frac{1}{p_i}$. $~~~~\blacksquare$

\textbf{Proof of Proposition 2.4.} By direct calculation, we have
\begin{align*}
\|f\|_{(L^{\vec{p}},L^{\vec{s}})}
&\sim\left\||B(\cdot,1)|^{\frac{1}{\alpha}-\frac{1}{n}
\sum^{n}_{i=1}\frac{1}{p_{i}}-\frac{1}{n}\sum^{n}_{i=1}\frac{1}{s_{i}}} \|f\chi_{B(\cdot,1)}\|_{L^{\vec{p}}} \right\|_{L^{\vec{s}}}\\
&\leq\sup_{r>0 }\left\||B(\cdot,r)|^{\frac{1}{\alpha}-\frac{1}{n}
\sum^{n}_{i=1}\frac{1}{p_{i}}-\frac{1}{n}\sum^{n}_{i=1}\frac{1}{s_{i}}} \|f\chi_{B(\cdot,r)}\|_{L^{\vec{p}}}\right\|_{L^{\vec{s}}}\\
&=\|f\|_{(L^{\vec{p}},L^{\vec{s}})^{\alpha}}.
\end{align*}
Therefore, $(L^{\vec{p}},L^{\vec{s}})^\alpha(\R^n)\subset(L^{\vec{p}},L^{\vec{s}})(\R^n)$ with $\|f\|_{(L^{\vec{p}},L^{\vec{s}})}\le\|f\|_{(L^{\vec{p}},L^{\vec{s}})^{\alpha}}$. Particularly, if $\vec{p}\le\vec{q}$, by H\"older's inequality,
$$|B(x,r)|^{\frac{1}{\alpha}-\frac{1}{n}\sum^{n}_{i=1}\frac{1}{p_{i}}-\frac{1}{n}\sum^{n}_{i=1}\frac{1}{s_{i}}} \|f\chi_{B(x,r)}\|_{L^{\vec{p}}}
\le|B(x,r)|^{\frac{1}{\alpha}-\frac{1}{n}\sum^{n}_{i=1}\frac{1}{q_{i}}-\frac{1}{n}\sum^{n}_{i=1}\frac{1}{s_{i}}} \|f\chi_{B(x,r)}\|_{L^{\vec{q}}} $$
Thus, $\|f\|_{(L^{\vec{p}},L^{\vec{s}})^{\alpha}}\le\|f\|_{(L^{\vec{q}},L^{\vec{s}})^{\alpha}}$ and $(L^{\vec{q}},L^{\vec{s}})^{\alpha}(\R^n)\subseteq(L^{\vec{p}},L^{\vec{s}})^{\alpha}(\R^n)$. $~~~~\blacksquare$

Before the proof of Proposition 2.5, the following two lemmas are necessary.

\textbf{Lemma 3.1.} Let $1\le\vec{p},\vec{s}\le\infty$ and $\frac{1}{n}\sum_{i=1}^n\frac{1}{s_i}\le\frac{1}{\alpha}\le\frac{1}{n}\sum_{i=1}^n\frac{1}{p_i}$. For any constant $\rho\in (0,\infty)$, we have
$$\left\|\|f\chi_{B(\cdot,r)}\|_{L^{\vec{p}}}\right\|_{L^{\vec{s}}}\sim \left\|\|f\chi_{B(\cdot,\rho r)}\|_{L^{\vec{p}}}\right\|_{L^{\vec{s}}}$$
where the positive equivalence constant are independent of $f$ and $t$.

\textbf{Proof.} Firstly, we prove the lemma holds when $\rho>1$. It is obvious that
$$\left\|\|f\chi_{B(\cdot,r)}\|_{L^{\vec{p}}}\right\|_{L^{\vec{s}}}\le\left\|\|f\chi_{B(\cdot,\rho r)}\|_{L^{\vec{p}}}\right\|_{L^{\vec{s}}}.$$
Next, we prove the reverse inequality. It is easy to find $N\in\mathbb{N}$ and $\{x_1,x_2,\cdots,x_N\}$, such that
$$B(0,\rho r)\subset\bigcup_{j=1}^{N}B(x_j,r),$$
where $N$ is independent of $r$ and $N\lesssim 1$. Therefore, we have
$$\|f\chi_{B(x,\rho r)}\|_{L^{\vec{p}}}\le\left\|f\sum_{j=1}^{N}\chi_{B(x+x_j,r)}\right\|_{L^{\vec{p}}} \le\sum_{j=1}^{N}\left\|f\chi_{B(x+x_j,r)}\right\|_{L^{\vec{p}}}$$
for any $x\in\R^n$. According to the translation invariance of the Lebesgue measure and $N\lesssim 1$, it follows that
$$\left\|\|f\chi_{B(\cdot,\rho r)}\|_{L^{\vec{p}}}\right\|_{L^{\vec{s}}} \le\sum_{j=1}^{N}\left\|\|f\chi_{B(\cdot+x_j,r)}\|_{L^{\vec{p}}}\right\|_{L^{\vec{s}}} \lesssim\left\|\|f\chi_{B(\cdot,r)}\|_{L^{\vec{p}}}\right\|_{L^{\vec{s}}}.$$

For the $\rho\in(0,1)$, we only need replace $r$ by $r/\rho$. The proof is completed. $~~~~\blacksquare$

\textbf{Remark 3.2.} If taking $r=1$, we have
$$\|f\|_{(L^{\vec{p}},L^{\vec{s}})}\sim\left\|\|f\chi_{B(\cdot,\rho)}\|_{L^{\vec{p}}}\right\|_{L^{\vec{s}}},~~\rho\in(0,\infty),$$
where the positive equivalence constants are independent of $f$.

The following result play an indispensable role in the proof of Proposition 2.5.

\textbf{Lemma 3.3.} Let $1\le\vec{p},\vec{s}\le\infty$ and $\frac{1}{n}\sum_{i=1}^n\frac{1}{s_i}\le\frac{1}{\alpha}\le\frac{1}{n}\sum_{i=1}^n\frac{1}{p_i}$. Then we have
$$\left\|\left\{\|f\chi_{Q_{r,k}}\|_{L^{\vec{p}}}\right\}_{k\in\mathbb{Z}^n}\right\|_{\ell^{\vec{s}}}\sim r^{-\sum_{i=1}^n\frac{1}{s_i}}\left\|\|f\chi_{B(\cdot,r)}\|_{L^{\vec{p}}}\right\|_{L^{\vec{s}}},$$
where the positive equivalence constants are independent of $f$ and $t$.

\textbf{Proof.} By the Lemma 3.1, we only need show that
$$\left\|\left\{\|f\chi_{Q_{r,k}}\|_{L^{\vec{p}}}\right\}_{k\in\mathbb{Z}^n}\right\|_{\ell^{\vec{s}}}\sim r^{-\sum_{i=1}^n\frac{1}{s}}\left\|\|f\chi_{B(\cdot,2\sqrt{n}r)}\|_{L^{\vec{p}}}\right\|_{L^{\vec{s}}}.$$
For any given $x\in\R^n$, we let
$$A_x:=\{k\in\mathbb{Z}:Q_{r,k}\cap B(x,2\sqrt{n}r)\neq\emptyset\}.$$
Then the cardinality of $A_x$ is finite and $x\in B(r k,4\sqrt{n}r)$ for any $k\in A_x$. Thus,
\begin{align*}
\|f\chi_{B(x,2\sqrt{n}r)}\|_{L^{\vec{p}}}&\le\left\|\sum_{k\in A_x}f\chi_{Q_{r,k}}\right\|_{L^{\vec{p}}}
\le\sum_{k\in A_x}\left\|f\chi_{Q_{r,k}}\right\|_{L^{\vec{p}}}\\
&\le\sum_{k\in \mathbb{Z}^n}\left\|f\chi_{Q_{r,k}}\right\|_{L^{\vec{p}}}\chi_{B(r k,4\sqrt{n}r)}(x).
\end{align*}
Taking $L^{\vec{s}}$-norm on $x$, we have
$$\left\|\|f\chi_{B(\cdot,2\sqrt{n}r)}\|_{L^{\vec{p}}}\right\|_{L^{\vec{s}}}\le\left\|\sum_{k\in \mathbb{Z}^n}\left\|f\chi_{Q_{r,k}}\right\|_{L^{\vec{p}}}\chi_{B(r k,4\sqrt{n}r)}\right\|_{L^{\vec{s}}}.$$
By the similar argument of Lemma 3.1, there exist $N\in\mathbb{N}$ and $\{k_1,k_2,\cdots,k_N\}$, such that
$$B(0,4\sqrt{n}r)\subset\bigcup_{j=1}^{N}Q_{k_j,r}$$
where $N$ is independent of $r$ and $N\sim 1$.
According to the translation invariance of the Lebesgue measure and $N\sim 1$, it follows that
\begin{align*}
\left\|\|f\chi_{B(\cdot,2\sqrt{n}r)}\|_{L^{\vec{p}}}\right\|_{L^{\vec{s}}}&\le\left\|\sum_{k\in \mathbb{Z}^n}\left\|f\chi_{Q_{r,k}}\right\|_{L^{\vec{p}}}\chi_{B(\rho k,4\sqrt{n}r)}\right\|_{L^{\vec{s}}}\\
&\le\left\|\sum_{j=1}^N\sum_{k\in \mathbb{Z}^n}\left\|f\chi_{Q_{r,k}}\right\|_{L^{\vec{p}}}\chi_{Q_{r,k_j+k}}\right\|_{L^{\vec{s}}}\\
&\lesssim r^{\sum_{i=1}^n\frac{1}{s_i}} \left\|\left\{\left\|f\chi_{Q_{r,k}}\right\|_{L^{\vec{p}}}\right\}_{k\in\mathbb{Z}^n}\right\|_{\ell^{\vec{s}}}.
\end{align*}
Indeed, the last inequality is obtained by the following fact that
\begin{align*}
&\quad\left(\int_{\R}\cdots\left(\int_{\R}\left| \sum_{k\in\mathbb{Z}^n}C_k\chi_{kr+(0,r]^n}(x)\right|^{s_1}dx_1\right)^{\frac{s_2}{s_1}}\cdots dx_n\right)^{\frac{1}{s_n}}\\
&=\left(\int_{\R}\cdots\left(\int_{\R} \left|\sum_{k\in\mathbb{Z}^n}C_k\prod_{i=1}^n\chi_{I_{k_i}}(x_i)\right|^{s_1}dx_1\right)^{\frac{s_2}{s_1}}\cdots dx_n\right)^{\frac{1}{s_n}}\\
&=\left(\sum_{k_n\in\mathbb{Z}}\int_{I_{k_n}}\cdots\left( \sum_{k_1\in\mathbb{Z}^n}\int_{I_{k_1}}\left|C_k\right|^{s_1}dx_1\right)^{\frac{s_2}{s_1}}\cdots dx_n\right)^{\frac{1}{s_n}}\\
&=r^{\sum_{i=1}^ns_i}\cdot\left(\sum_{k_n\in\mathbb{Z}}\cdots\left( \sum_{k_1\in\mathbb{Z}^n}\left|C_k\right|^{s_1}dx_1\right)^{\frac{s_2}{s_1}}\cdots dx_n\right)^{\frac{1}{s_n}},
\end{align*}
where $C_k=\left\|f\chi_{Q_{r,k}}\right\|_{L^{\vec{p}}}$ and $I_{k_i}=rk_i+(0,r]$. Thus, we prove that
$$r^{-\sum_{i=1}^n\frac{1}{s_i}}\left\|\|f\chi_{B(\cdot,2\sqrt{n}r)}\|_{L^{\vec{p}}}\right\|_{L^{\vec{s}}} \lesssim \left\|\left\{\left\|f\chi_{Q_{r,k}}\right\|_{L^{\vec{p}}}\right\}_{k\in\mathbb{Z}^n}\right\|_{\ell^{\vec{s}}}.$$
For the opposite inequality, it is obvious that
$$r^{\sum_{i=1}^n\frac{1}{s_i}}\left\|\left\{\left\|f\chi_{Q_{r,k}}\right\|_{L^{\vec{p}}}\right\}_{k\in\mathbb{Z}^n}\right\|_{\ell^{\vec{s}}}
=\left\|\sum_{k\in\mathbb{Z}^n}\left\|f\chi_{Q_{r,k}}\right\|_{L^{\vec{p}}}\chi_{Q_{r,k}}\right\|_{L^{\vec{s}}}.$$
By $Q_{r,k}\subset B(x,2\sqrt{n}r)$ for $x\in Q_{r,k}$, we have
$$r^{\sum_{i=1}^n\frac{1}{s_i}}\left\|\left\{\left\|f\chi_{Q_{r,k}}\right\|_{L^{\vec{p}}}\right\}_{k\in\mathbb{Z}^n}\right\|_{\ell^{\vec{s}}} \le\left\|\left\|f\chi_{B(\cdot,2\sqrt{n}r)}\right\|_{L^{\vec{p}}}\right\|_{L^{\vec{s}}}.$$
The proof is completed. $~~~~\blacksquare$

By Lemma 3.3, the proof of Proposition 2.5 is easy.

\textbf{Proof of Proposition 2.5.} According to the Lemma 3.3, we obtain that
$$\left\|\left\{\|f\chi_{Q_{1,k}}\|_{L^{\vec{p}}}\right\}_{k\in\mathbb{Z}^n}\right\|_{\ell^{\vec{s}}}\sim \left\|\|f\chi_{B(\cdot,1)}\|_{L^{\vec{p}}}\right\|_{L^{\vec{s}}}$$
and
$$r^{\frac{n}{\alpha}-\sum_{i=1}^{n}\frac{1}{p_i}} \left\|\left\{\|f\chi_{Q_{r,k}}\|_{L^{\vec{p}}}\right\}_{k\in\mathbb{Z}^n}\right\|_{L^{\vec{s}}}\sim \left\||B(\cdot,r)|^{\frac{1}{\alpha}-\frac{1}{n}\sum_{i=1}^{n}\frac{1}{p_i}-\frac{1}{n}\sum_{i=1}^n\frac{1}{s_i}} \|f\chi_{B(\cdot,r)}\|_{L^{\vec{p}}}\right\|_{L^{\vec{s}}}.$$
Thus, we prove proposition 2.5.  $~~~~\blacksquare$

\section{The proof of Theorem 2.7}\label{sec3}
\par
In this section, we will prove Theorem 2.7, whose the ideal comes from \cite{33}. Before that, the dual of mixed-norm amalgam spaces $(L^{\vec{p}},L^{\vec{s}})(\R^n)$ will given as follows.

\textbf{Lemma 4.1.} (\romannumeral1) Let $1\le\vec{p},\vec{s}\le\infty$. For $r\in(0,\infty)$, we have
$$\|fg\|_1\le{_r\|f\|_{\vec{p},\vec{s}}}\cdot{_r\|g\|_{\vec{p}',\vec{s}\,'}},~~f,g\in L_{loc}^{1}(\R^n).\eqno{(4.1)}$$
(\romannumeral2) Let $1\le\vec{p},\vec{s}<\infty$. The dual of mixed-norm amalgam spaces $(L^{\vec{p}},\ell^{\vec{s}})(\R^n)$ is $(L^{\vec{p}'},\ell^{\vec{s}\,'})(\R^n)$.

\textbf{Proof.} For $0<r<\infty$, by H\"older's inequality, we have
$$\|fg\|_1\le{_r\|f\|_{\vec{p},\vec{s}}}\cdot{_r\|g\|_{\vec{p}',\vec{s}\,'}},~~f,g\in L_{loc}^{1}(\R^n).$$
According to Theorem 2 of \cite{5} and Theorem 1a) of Section 3 in \cite{11}, it immediate to deduce that the dual of $(L^{\vec{p}},\ell^{s_1})(\R^n)$ is $(L^{\vec{p}'},\ell^{s'_1})(\R^n)$. If the dual of $(L^{\vec{p}},\ell^{\bar{s}})(\R^n)$ is $(L^{\vec{p}'},\ell^{\bar{s}'})(\R^n)$ with $\bar{s}=(s_1,s_2,\cdots,s_{n-1})$, using Theorem 2 of \cite{5},
$$(L^{\vec{p}},\ell^{\vec{s}})^{*}=\left(\prod_{k_n\in\mathbb{Z}}(L^{\vec{p}},\ell^{\bar{s}}),\ell^{s_n}\right)^{*} =\left(\prod_{k_n\in\mathbb{Z}}(L^{\vec{p}},\ell^{\bar{s}})^{*},\left(\ell^{s_n}\right)^{*}\right) =\left(\prod_{k_n\in\mathbb{Z}}(L^{\vec{p}'},\ell^{\bar{s}'}),\ell^{s'_n}\right)=(L^{\vec{p}'},\ell^{\vec{s}\,'}).$$
Hence, $(L^{\vec{p}'},\ell^{\vec{s}\,'})(\R^{n})$ is isometrically isomorphic to the dual of $(L^{\vec{p}},\ell^{\vec{s}})(\R^{n})$. There is an unique element $\phi(T)$ of $(L^{\vec{p}},\ell^{\vec{s}})(\R^{n})$ such that
$$T(f)=\int_{\R^n}f(x)\phi(T)(x)dx,~~f\in(L^{\vec{p}},\ell^{\vec{s}})(\R^{n})$$
and furthermore
$$\|\phi(T)\|_{\vec{p}',\vec{s}\,'}=\|T\|,\eqno{(4.2)}$$
where $\|T\|:=\sup\left\{|T(f)|:f\in L_{loc}^{1}(\R^n) \text{ and }\|f\|_{\vec{p},\vec{s}}\le 1\right\}$. $~~~~\blacksquare$

Now, we discuss the properties of the dilation operator $St_{r}^{(\alpha)}:f\mapsto r^{-\frac{n}{\alpha}}f(r^{-1}\cdot)$ for $0<\alpha<\infty$ and $0<r<\infty$. By direct computation, we have the following properties.

\textbf{Proposition 4.2.} Let $f\in L^1_{loc}(\R^n)$, $0<\alpha<\infty$, and $0<r<\infty$.\\
(\romannumeral1) $St_{r}^{(\alpha)}$ maps $L^1_{loc}(\R^n)$ into itself.\\
(\romannumeral2) $f=St_{1}^{(\alpha)}(f)$.\\
(\romannumeral3) $St_{r_1}^{(\alpha)}\circ St_{r_2}^{(\alpha)}=St_{r_2}^{(\alpha)}St_{r_2}^{(\alpha)}=St_{r_1r_2}^{(\alpha)}$.\\
(\romannumeral4) $\sup_{r>0}\|St_{r}^{(\alpha)}(f)\|_{\vec{p},\vec{s}}=\|f\|_{\vec{p},\vec{s},\alpha}$, where $1\le\vec{p},\vec{s}<\infty$ and $\frac{1}{n}\sum_{i=1}^n\frac{1}{s_i}\le\frac{1}{\alpha}\le\frac{1}{n}\sum_{i=1}^n\frac{1}{p_i}$.

By Proposition 4.2 and Definition 2.6, the following result can be obtained.

\textbf{Proposition 4.3.} Let $1\le\vec{p},\vec{s}\le\infty$, and $\frac{1}{n}\sum_{i=1}^n\frac{1}{s_i}\leq\frac{1}{\alpha}\leq\frac{1}{n}\sum^{n}_{i=1}\frac{1}{p_{i}}$. $(L^{\vec{p}\,'},\ell^{\vec{s}\,'})(\R^n)$ is a dense subspace of $\mathcal{H}(\vec{p}\,',\vec{s}\,',\alpha')$.

\textbf{Proof.} First we verify that $(L^{\vec{p}\,'},\ell^{\vec{s}\,'})(\R^n)$ is continuously embedded into $\mathcal{H}(\vec{p}\,',s',\alpha')$. For any $f\in(L^{\vec{p}\,'},\ell^{\vec{s}\,'})(\R^n)$, we have
$$f=\|f\|_{\vec{p}\,',\vec{s}\,'}St_1^{\alpha}(\|f\|^{-1}_{\vec{p}\,',\vec{s}\,'}f)\eqno{(4.3)}$$
and
$$\left\|\|f\|^{-1}_{\vec{p}\,',\vec{s}\,'}f\right\|_{\vec{p}\,',\vec{s}\,'}=1.$$
Thus, $f\in\mathcal{H}(\vec{p}\,',\vec{s}\,',\alpha')$ and satisfies
$$\|f\|_{\mathcal{H}(\vec{p}\,',\vec{s}\,',\alpha')}\le\|f\|_{\vec{p}\,',\vec{s}\,'}.$$

Let us show the denseness of $(L^{\vec{p}\,'},\ell^{\vec{s}\,'})(\R^n)$ in $\mathcal{H}(\vec{p}\,',\vec{s}\,',\alpha')$. It is clear that if $\{(c_j,r_j,f_j)\}_{j\ge 1}$ is a block decomposition of $f\in\mathcal{H}(\vec{p}\,',\vec{s}\,',\alpha')$, then
$$\left\{\sum_{j=1}^Jc_j St_{r_j}^{(\alpha')}(f_j)\right\}_{J\ge 1}\subset(L^{\vec{p}\,'},\ell^{\vec{s}\,'})(\R^n)$$
and
$$\left\|f-\sum_{j=1}^{J}c_jSt_{r_j}^{(\alpha')}(f_j)\right\|_{\mathcal{H}(\vec{p}\,',\vec{s}\,',\alpha')} =\left\|\sum_{j=J+1}^{\infty}c_jSt_{r_j}^{(\alpha')}(f_j)\right\|_{\mathcal{H}(\vec{p}\,',\vec{s}\,',\alpha')} \le\sum_{j=J+1}^{\infty}|c_j|\rightarrow 0$$
with $J\rightarrow\infty$. Thus, $(L^{\vec{p}\,'},\ell^{\vec{s}\,'})(\R^n)$ is a dense subspace of $\mathcal{H}(\vec{p}\,',\vec{s}\,',\alpha')$. $~~~~\blacksquare$

Now, let us to prove the main theorem in this section.

\textbf{The proof of Theorem 2.7.} Let us to prove (\romannumeral1). Let $\{(c_j,r_j,f_j)\}_{j\ge 1}$ be block decomposition of $f$. By Proposition 4.1 and (4.1), we have for any $j\ge 1$
\begin{align*}
\left|\int_{\R^n} St_{r_j}^{(\alpha')}(f_j)(x)g(x)dx\right|&=\left|\int_{\R^n} St_{r^{-1}_j}^{(\alpha)}(g)(x)f_j(x)dx\right|\\
&\le\int_{\R^n}\left|St_{r^{-1}_j}^{(\alpha)}(g)(x)f_j(x)\right|dx\\
&\le\|f_j\|_{\vec{p}',\vec{s}\,'}\left\|St_{r^{-1}_j}^{(\alpha)}(g)\right\|_{\vec{p},\vec{s}}\\
&\le\left\|St_{r^{-1}_j}^{(\alpha)}(g)\right\|_{\vec{p},\vec{s}}\le\|g\|_{\vec{p},\vec{s},\alpha}.
\end{align*}
Therefore we have
$$\sum_{j\ge 1}\int_{\R^n}\left|c_jSt_{r_j}^{(\alpha')}(f_j)(x)g(x)\right|dx\le\|g\|_{\vec{p},\vec{s},\alpha}\sum_{j\ge 1}|c_j|.$$
This implies that $fg=g\sum_{j\le 1}c_jSt_{r_j}^{(\alpha')}(f_j)$ belong to $L^1(\R^n)$ and
$$\left|\int_{\R^n}f(x)g(x)dx\right|\le\int_{\R^n}|f(x)g(x)|dx\le\|g\|_{\vec{p},\vec{s},\alpha}\sum_{j\ge 1}|c_j|.$$
Taking the infimum with respect to all block decompositions of $f$, we get
$$\left|\int_{\R^n}f(x)g(x)dx\right|\le\int_{\R^n}|f(x)g(x)|dx\le\|g\|_{\vec{p},\vec{s},\alpha}\|f\|_{\vec{p}',\vec{s}\,',\alpha'}.$$

Now, Let us prove (\romannumeral2). By the (\romannumeral1), we have
$$T_g\in\mathcal{H}(\vec{p},\vec{s},\alpha)^{*}.$$
For any $a_1,a_2\in R,~g_1,g_2\in(L^{\vec{p}},\ell^{\vec{s}})^{\alpha}(\R^n)$
$$T(a_1g_1+a_2g_2)=a_1T_{g_1}+a_2T_{g_2}$$
and
$$\|T_g\|=\sup_{\|f\|_{\mathcal{H}(\vec{p},\vec{s},\alpha)}\le 1}|T_g(f)|\le\|g\|_{\vec{p},\vec{s},\alpha},$$
that is, $T$ is linear and bounded mapping from $(L^{\vec{p}},\ell^{\vec{s}})^{\alpha}(\R^n)$ into $\mathcal{H}{(\vec{p}',\vec{s}\,',\alpha')}^{*}$ satisfying $\|T\|\le 1$. For any $g_1,g_2\in(L^{\vec{p}},\ell^{\vec{s}})^{\alpha}(\R^n)\subset(L^{\vec{p}},\ell^{\vec{s}})(\R^n)$, if $T_{g_1}=T_{g_2}$, then for any $f\in(L^{\vec{p}\,'},\ell^{\vec{s}\,'})(\R^n)\subset\mathcal{H}(\vec{p}\,',\vec{s}\,',\alpha')$, we have
$$T_{g_1}(f)=T_{g_2}(f).$$
Thus, $g_1=g_2$, that is , $T$ is injective.

Now, we will prove that $T$ is a surjection and $\|g\|_{\vec{p},\vec{s},\alpha}\le\|T_g\|~(\text{or }\|T\|\ge 1)$. Let $T^{*}$ be an element of $\mathcal{H}(\vec{p}\,',\vec{s}\,',\alpha')^{*}$. From Proposition 4.3, it follows that the restriction $T^{*}_0$ of $T^{*}$ to $(L^{\vec{p}\,'},\ell^{\vec{s}\,'})(\R^n)$ belong to $\mathcal{H}(\vec{p}\,',\vec{s}\,',\alpha')^{*}$. Furthermore, we have $$\frac{1}{n}\sum_{j=1}^{n}\frac{1}{p_i'}\leq\frac{1}{\alpha'}\le\frac{1}{n}\sum_{i=1}^n\frac{1}{s'_i}.$$
There is an element $g$ of $(L^{\vec{p}},\ell^{\vec{s}})(\R^n)$ such that for any $f\in(L^{\vec{p}\,'},\ell^{\vec{s}\,'})(\R^n)$
$$T^{*}(f)=T_0^{*}(f)=\int_{\R^n}f(x)g(x)dx.\eqno{(4.4)}$$
Hence, for $f\in(L^{\vec{p}\,'},\ell^{\vec{s}\,'})(\R^n)$ and $\rho>0$ we have
$$\int_{\R^n}St^{(\alpha)}_{r}(g)(x)f(x)dx=\int_{\R^n}g(x)St^{(\alpha')}_{r^{-1}}(f)(x)dx=T^{*}\left[St^{(\alpha')}_{r^{-1}}(f)\right].$$
and $St^{(\alpha')}_{r^{-1}}(f)\in\mathcal{H}(\vec{p}\,',\vec{s}\,',\alpha')$. By the the assumption $T^{*}\in\mathcal{H}(\vec{p}\,',\vec{s}\,',\alpha')^{*}$, we have
$$\left|\int_{\R^n}St^{(\alpha)}_{r}(g)(x)f(x)dx\right| \le\|T^{*}\|\cdot\|St^{(\alpha')}_{\rho^{-1}}(f)\|_{\mathcal{H}(\vec{p}\,',\vec{s}\,',\alpha')} \le\|T^{*}\|\cdot\|f\|_{\vec{p}\,',\vec{s}\,'}.$$
Due to (4.2), it follows that
$$\|St^{(\alpha)}_{r}(g)\|_{\vec{p},\vec{s}}\le\|T^{*}\|.$$
Therefore, for any $g\in(L^{\vec{p}},\ell^{\vec{s}})(\R^n)$, by Proposition 4.2,
$$\|g\|_{\vec{p},\vec{s},\alpha}\le\|T^{*}\|.$$
According to (4.4) and Proposition 4.3, we get
$$T^{*}(f)=\int_{\R^n}f(x)g(x)dx,~~f\in\mathcal{H}(\vec{p}\,',\vec{s}\,',\alpha').$$
Thus, $T$ is a surjection and $\|g\|_{\vec{p},\vec{s},\alpha}\le\|T\|$.

\section{The proof of Theorem 2.8}\label{sec3}
\par
In this section, we will prove the conclusions of Theorem 2.8.

\textbf{Proof of Theorem 2.8.} By Remark 2.9, we only need to prove the boundedness of $I_\gamma$ on mixed-norm amalgam spaces if $\gamma=\frac{n}{\alpha}-\frac{n}{\beta}$. Let $f\in(L^{\vec{q}},L^{\vec{s}})^{\alpha}(\R^n)$, $B=B(y,r)$, and
$$f=f_{1}+f_{2}=f\chi_{2B}+f\chi_{(2B)^{c}}.$$
where $\chi_{2B}$ is the characteristic function of $2B$. By the linearity of the fractional integral operator $I_\gamma$, one can write
\begin{align*}
&\quad|B(y,r)|^{\frac{1}{\beta}-\frac{1}{n}\sum_{i=1}^{n}\frac{1}{q_{i}}-\frac{1}{n}\sum_{i=1}^n\frac{1}{s_i}}
\|I_{\alpha}(f)\chi_{B(y,r)}\|_{L^{\vec{q}}}\\
&=|B(y,r)|^{\frac{1}{\beta}-\frac{1}{n}\sum_{i=1}^{n}\frac{1}{q_{i}}-\frac{1}{n}\sum_{i=1}^n\frac{1}{s_i}}
\|I_{\alpha}(f_1)\chi_{B(y,r)}\|_{L^{\vec{q}}}\\
&+|B(y,r)|^{\frac{1}{\beta}-\frac{1}{n}\sum_{i=1}^{n}\frac{1}{q_{i}}-\frac{1}{n}\sum_{i=1}^n\frac{1}{s_i}}
\|I_{\alpha}(f_2)\chi_{B(y,r)}\|_{L^{\vec{q}}}\\
&:=\text{I}(y,r)+\text{II}(y,r)
\end{align*}

Below, we will give the estimates of $I(y,r)$ and $II(y,r)$, respectively. By the$(L^{\vec{p}},L^{\vec{q}})$-boundedness of $I_\gamma$ (see Lemma 1.1),
\begin{align*}
\text{I}(y,r)&=|B(y,r)|^{\frac{1}{\beta}-\frac{1}{n}\sum_{i=1}^{n}\frac{1}{q_{i}}-\frac{1}{n}\sum_{i=1}^n\frac{1}{s_i}}
\|I_{\alpha}(f_1)\chi_{B(y,r)}\|_{L^{\vec{q}}}\\
&\le|B(y,r)|^{\frac{1}{\beta}-\frac{1}{n}\sum_{i=1}^{n}\frac{1}{q_{i}}-\frac{1}{n}\sum_{i=1}^n\frac{1}{s_i}} \|f\chi_{2B(y,r)}\|_{L^{\vec{p}}}\\
&\sim|2B(y,r)|^{\frac{1}{\alpha}-\frac{1}{n}\sum_{i=1}^{n}\frac{1}{p_{i}}-\frac{1}{n}\sum_{i=1}^n\frac{1}{s_i}}
\|f\chi_{2B(y,r)}\|_{L^{\vec{p}}}.
\end{align*}
Thus,
$$\sup_{r>0} \|\text{I}(\cdot,r)\|_{L^{\vec{s}}}\lesssim\|f\|_{(L^{\vec{p}},L^{\vec{s}})^{\alpha}}.\eqno{(5.1)}$$

Let us now turn to the estimate of $\text{II}(y,r)$. First, it is clear that when $x\in B(y,r)$ and $z\in (2B)^{c}$, we get $|x-z|\sim|y-z|$. Then we decompose $\R^n$ into a geometrically increasing sequence of concentric balls and obtain the following pointwise estimate:
\begin{align*}
I_{\gamma}(f_{2})(x)&=\int_{R^{n}}\frac{|f_{2}(z)|}{|x-z|^{n-\gamma}}dz\\
&=\int_{(2B)^{c}}\frac{|f(z)|}{|x-z|^{n-\gamma}}dz\qquad\qquad\\
&\sim\sum\limits_{j=1}^{\infty}\int_{2^{j+1}B\setminus2^{j}B}\frac{|f(z)|}{|x-z|^{n-\gamma}}dz\qquad\\
&\lesssim\sum_{j=1}^{\infty}\frac{1}{|2^{j+1}B|^{1-\frac{\gamma}{n}}}\int_{2^{j+1}B}|f(z)|dz.
\end{align*}
Combining
$$I_{\gamma}(f_{2})(x)\lesssim\sum_{j=1}^{\infty}\frac{1}{|2^{j+1}B(y,r)|^{1-\frac{\gamma}{n}}}\int_{2^{j+1}B(y,r)}|f(z)|dz\eqno{(5.2)}$$
and H\"older's inequality, we obtain
\begin{eqnarray*}
&&\text{II}=|B(y,r)|^{\frac{1}{\beta}-\frac{1}{n}\sum_{i=1}^{n}\frac{1}{q_{i}}-\frac{1}{n}\sum_{i=1}^n\frac{1}{s_i}}
\left\|\chi_{B(y,r)}\int_{R^{n}}\frac{|f_{2}(z)|}{|\cdot-z|^{n-\gamma}}dz\right\|_{L^{\vec{q}}}\\
&&\quad\lesssim|B(y,r)|^{\frac{1}{\beta}-\frac{1}{n}\sum\limits_{i=1}^{n}\frac{1}{q_{i}}-\frac{1}{n}\sum_{i=1}^n\frac{1}{s_i}}
\sum\limits_{j=1}^{\infty}\frac{1}{|2^{j+1}B(y,r)|^{1-\frac{\gamma}{n}}}\int_{2^{j+1}B(y,r)}|f(z)|dz
|B(y,r)|^{\frac{1}{n}\sum\limits_{i=1}^{n}\frac{1}{q_{i}}}\\
&&\quad\lesssim\sum\limits_{j=1}^{\infty}
|B(y,r)|^{\frac{1}{\beta}-\frac{1}{n}\sum_{i=1}^n\frac{1}{s_i}}\frac{1}{|2^{j+1}B(y,r)|^{1-\frac{\gamma}{n}}}
|2^{j+1}B(y,r)|^{\frac{1}{n}\sum_{i=1}^{n}\frac{1}{p'_{i}}}\|f\chi_{2^{j+1}B(y,r)}\|_{L^{\vec{p}}}\\
&&\quad\sim\sum\limits_{j=1}^{\infty}
2^{-j(\frac{1}{\beta}-\frac{1}{n}\sum_{i=1}^n\frac{1}{s_i})}
|2^{j+1}B(y,r)|^{\frac{1}{\beta}-\frac{1}{n}\sum_{i=1}^n\frac{1}{s_i}-\frac{1}{n}\sum_{i=1}^{n}\frac{1}{p_{i}}
+1-1+\frac{\gamma}{n}}
\|f\chi_{2^{j+1}B(y,r)}\|_{L^{\vec{p}}}\\
&&\quad=\sum_{j=1}^{\infty} 2^{-j(\frac{1}{\beta}-\frac{1}{n}\sum_{i=1}^n\frac{1}{s_i})}
|2^{j+1}B(y,r)|^{\frac{1}{\alpha}-\frac{1}{n}\sum_{i=1}^n\frac{1}{s_i}-\frac{1}{n}\sum_{i=1}^{n}\frac{1}{p_{i}}}
\|f\chi_{2^{j+1}B(y,r)}\|_{L^{\vec{p}}}
\end{eqnarray*}
By $\frac{1}{\beta}-\frac{1}{n}\sum_{i=1}^n\frac{1}{s_i}>0$,
$$\sum_{j=1}^{\infty}2^{-j(\frac{1}{\beta}-\frac{1}{n}\sum_{i=1}^n\frac{1}{s_i})}\sim 1.$$
Thus,
$$\sup_{r>0}\|\text{II}(\cdot,r)\|_{L^{\vec{s}}}\lesssim\|f\|_{(L^{\vec{p}},L^{\vec{s}})^{\alpha}}.\eqno{(5.3)}$$
Therefore, using (5.1) and (5.3),
\begin{align*}
\|I_{\alpha}f\|_{(L^{\vec{p}},L^{\vec{s}})^{\beta}} &=\sup_{r>0}\left\||B(\cdot,r)|^{\frac{1}{\beta}-\frac{1}{n}\sum_{i=1}^{n}\frac{1}{q_{i}}-\frac{1}{n}\sum_{i=1}^n\frac{1}{s_i}}
\|I_{\alpha}(f)\chi_{B(\cdot,r)}\|_{L^{\vec{q}}}\right\|_{L^{\vec{s}}}\\
&\le\sup_{r>0}\|\text{I}(\cdot,r)\|_{L^{\vec{s}}}+\sup_{r>0}\|\text{II}(\cdot,r)\|_{L^{\vec{s}}}\\
&\lesssim\|f\|_{(L^{\vec{p}},L^{\vec{s}})^{\alpha}}.
\end{align*}
The proof is completed. $~~~~\blacksquare$

Let $0<\gamma<n$. The related fractional maximal function is defined as
$$M_{\gamma}f(x):=\sup_{B\ni x}\frac{1}{|B|^{1-\frac{\gamma}{n}}}\int_{B}|f(y)|dy,$$
where the supremum is taken over all cube $B\subset\mathbb{R}^n$ containing $x$. It is well-know that
$$|M_{\gamma}f(x)|\lesssim I_{\gamma}(|f|)(x).\eqno{(5.4)}$$
An immediate application of the above inequality (5.4) is the following strong-type for the operators $M_{\gamma}$.

\textbf{Corollary 5.1.} Let $0<\gamma<n$, $1<\vec{p},\vec{q}<\infty$, $1<\vec{s}\le\infty$, $\frac{1}{n}\sum_{i=1}^{n}\frac{1}{s_i}\le\frac{1}{\alpha}\le\frac{1}{n}\sum_{i=1}^n\frac{1}{p_i}$, and $\frac{1}{n}\sum_{i=1}^{n}\frac{1}{s_i}\le\frac{1}{\beta}\le\frac{1}{n}\sum_{i=1}^n\frac{1}{q_i}$. Assume that $\gamma=\sum_{i=1}^n\frac{1}{p_i}-\sum_{i=1}^{n}\frac{1}{q_i}=\frac{n}{\alpha}-\frac{n}{\beta}$. Then the fractional integral operators $M_\gamma$ are bounded from $(L^{\vec{p}},L^{\vec{s}})^{\alpha}(\R^n)$ to $(L^{\vec{q}},L^{\vec{s}})^{\beta}(\R^n)$.

Before the next corollary, let us recall generalized fractional integral operators.

Suppose that $\mathcal{L}$ are linear operators which generate an analytic semigroup $\{e^{-t\mathcal{L}}\}_{t>0}$ on
$L^2(\mathbb{R}^n)$ with a kernel $p_t(x,y)$ satisfying
$$|p_t(x,y)|\le\frac{C_1}{t^{n/2}}e^{-C_2\frac{|x-y|^2}{t}}~~x,y\in \mathbb{R}^n,$$
where $C_1,C_2>0$ are independent of $x,~y$ and $t$.

For any $0<\gamma<n$, the generalized fractional integral operators $\mathcal{L}^{-\gamma/2}$ associated with the operator $\mathcal{L}$ is defined by
$$\mathcal{L}^{-\gamma/2}f(x)=\frac{1}{\Gamma(\gamma/2)}\int_0^\infty e^{-t\mathcal{L}}(f)(x)\frac{dt}{t^{-\gamma/2+1}}.$$

Note that if $\mathcal{L}=-\Delta$ is the Laplacian on $\mathbb{R}^n$, then $\mathcal{L}^{-\gamma/2}$ is the classical fractional integral operators $I_{\gamma}$. See, for example, Chapter 5 of \cite{23}. By Gaussian upper bound of kernel $p_t(x,y)$, it is easy to check that for all $x\in\mathbb{R}^n$,
$$|\mathcal{L}^{-\gamma/2}f(x)|\le CI_{\gamma}(|f|)(x).$$
(see \cite{24}). In fact, if we denote the the kernel of $\mathcal{L}^{-\gamma/2}$ by $K_{\gamma}(x,y)$, it is easy to obtain that
\begin{align*}
\mathcal{L}^{-\gamma/2}f(x)&=\frac{1}{\Gamma(\gamma/2)}\int_0^\infty e^{-t\mathcal{L}}(f)(x)\frac{dt}{t^{-\gamma/2+1}}\\
&=\frac{1}{\Gamma(\gamma/2)}\int_0^\infty\int_{\mathbb{R}^n}p_t(x,y)f(y)dy\frac{dt}{t^{-\gamma/2+1}}\\
&=\int_{\mathbb{R}^n}\frac{1}{\Gamma(\gamma/2)}\int_0^\infty p_t(x,y)\frac{dt}{t^{-\gamma/2+1}}\cdot f(y)dy\\
&=\int_{\mathbb{R}^n}K_{\gamma}(x,y)\cdot f(y)dy.
\end{align*}
Hence, by Gaussian upper bound,
\begin{align*}
|K_{\gamma}(x,y)|&=\left|\frac{1}{\Gamma(\gamma/2)}\int_0^\infty p_t(x,y)\frac{dt}{t^{-\gamma/2+1}}\right|\\
&\le \frac{1}{\Gamma(\gamma/2)}\int_0^\infty |p_t(x,y)|\frac{dt}{t^{-\gamma/2+1}}\\
&\le C\int_0^\infty e^{-C_2\frac{|x-y|^2}{t}}\frac{dt}{t^{n/2-\gamma/2+1}}\\
&\le C\cdot\frac{1}{|x-y|^{n-\gamma}}.
\end{align*}
Taking into account this pointwise inequality, as a consequence of Theorem 2.3, we have the following result.

\textbf{Corollary 5.2} Let $0<\gamma<n$, $1<\vec{p},\vec{q}<\infty$, $1<\vec{s}\le\infty$, $\frac{1}{n}\sum_{i=1}^{n}\frac{1}{s_i}\le\frac{1}{\alpha}\le\frac{1}{n}\sum_{i=1}^n\frac{1}{p_i}$, and $\frac{1}{n}\sum_{i=1}^{n}\frac{1}{s_i}\le\frac{1}{\beta}\le\frac{1}{n}\sum_{i=1}^n\frac{1}{q_i}$. Assume that $\gamma=\sum_{i=1}^n\frac{1}{p_i}-\sum_{i=1}^{n}\frac{1}{q_i}=\frac{n}{\alpha}-\frac{n}{\beta}$. Then the generalized fractional integral operators $\mathcal{L}^{\gamma/2}$ are bounded from $(L^{\vec{p}},L^{\vec{s}})^{\alpha}(\R^n)$ to $(L^{\vec{q}},L^{\vec{s}})^{\beta}(\R^n)$.

\section{Proof of Theorem 2.10}\label{sec4}
\par
To prove Theorem 2.10 in this section, we need the following lemmas about $BMO(\R^n)$ function.

\textbf{Lemma 6.1} Let $b$ be a function in $BMO(\R^n)$.\\
(i) For any ball $B$ in $\R^n$ and for any positive integer $j\in \mathbb{Z}^{+}$,
$$|b_{2^{j+1}B}-b_B|\le Cj\|b\|_{*}.$$
(ii) Let $1<\vec{p}<\infty$. There exist positive constants $C_1\le C_2$ such that for all $b\in BMO(\R^n)$,
$$C_1\|b\|_{*}\le \sup_{B\subset\mathbb{R}^n}\frac{\|b-b_B\|_{L^{\vec{p}}(\R^n)}}{\|\chi_{B}\|_{L^{\vec{p}}(\R^n)}}\le C_2\|b\|_{*}.$$

\textbf{Proof.} The proof of (i) is so easy that we omit.  By Lemma 3.5 of \cite{22}, the $Mf$ is bounded on $L^{\vec{p}}(\mathbb{R}^n)$ with $1<\vec{p}=(p_1,p_2,\cdots,p_n)<\infty$. According to the dual theorem of  Theorem 1.a of \cite{11}, the associate spaces of $L^{\vec{p}}(\R^n)$ is $L^{\vec{p}\,'}(\mathbb{R}^n)$. Finally, by Theorem 1.1 of \cite{21},  the proof of (ii) can be proved.

Now, let us show the proof of Theorem 2.4.

\textbf{Proof of Theorem 2.10.} Let $f\in(\L^{\vec{p}},L^{\vec{s}})^{\alpha}(\R^n),~B=B(y,r)$ and
$$f=f_{1}+f_{2}=f\chi_{2B}+f\chi_{(2B)^{c}}.$$
By the linearity of the commutator operators $[b,I_{\gamma}]$, we write
\begin{align*}
&\quad|B(y,r)|^{\frac{1}{\beta}-\frac{1}{n}\sum_{i=1}^n\frac{1}{s_i}-\frac{1}{n}\sum_{i=1}^{n}\frac{1}{q_{i}}}
\|[b,I_{\gamma}](f)\chi_{B(y,r)}\|_{L^{\vec{q}}}\\
&\leq |B(y,r)|^{\frac{1}{\beta}-\frac{1}{n}\sum_{i=1}^n\frac{1}{s_i}-\frac{1}{n}\sum_{i=1}^{n}\frac{1}{q_{i}}}
\|[b,I_{\gamma}](f_{1})\chi_{B(y,r)}\|_{L^{\vec{q}}}\\
&+|B(y,r)|^{\frac{1}{\alpha}-\frac{1}{n}\sum_{i=1}^n\frac{1}{s_i}-\frac{1}{n}\sum_{i=1}^{n}\frac{1}{q_{i}}}
\|[b,I_{\gamma}](f_{2})\chi_{B(y,r)}\|_{L^{\vec{q}}}\\
&:=\text{I}(y,r)+\text{II}(y,r).
\end{align*}
By Lemma 1.2 and observe that $\frac{1}{\beta}-\frac{1}{n}\sum_{i=1}^{n}\frac{1}{q_i}=\frac{1}{\alpha}-\frac{1}{n}\sum_{i=1}^{n}\frac{1}{p_i}$,
\begin{align*}
\text{I}(y,r) &=|B(y,r)|^{\frac{1}{\beta}-\frac{1}{n}\sum_{i=1}^n\frac{1}{s_i}-\frac{1}{n}\sum_{i=1}^{n}\frac{1}{q_{i}}}
\|[b,I_{\gamma}](f_{1})\chi_{B(y,r)}\|_{L^{\vec{q}}}\\
&\lesssim|2B(y,r)|^{\frac{1}{\beta}-\frac{1}{n}\sum_{i=1}^n\frac{1}{s_i}-\frac{1}{n}\sum_{i=1}^{n}\frac{1}{q_{i}}} \|f\chi_{2B(y,r)}\|_{L^{\vec{p}}}\\
&=|2B(y,r)|^{\frac{1}{\alpha}-\frac{1}{n}\sum_{i=1}^n\frac{1}{s_i}-\frac{1}{n}\sum_{i=1}^{n}\frac{1}{p_{i}}}
\|f \chi_{2B(y,r)}\|_{L^{\vec{p}}}.
\end{align*}
Thus,
$$\sup_{r>0}\|I(y,r)\|_{L^{\vec{s}}}\lesssim\|f\|_{(L^{\vec{p}},L^{\vec{s}})^{\alpha}}.\eqno{(6.1)}$$

Now, let us turn to the estimate of $\text{II}(y,r)$. By the definition of $[b,I_{\gamma}]$, we have
$$|[b,I_{\gamma}](f_{2})(x)|\le|b(x)-b_{B(y,r)}|\cdot|I_{\gamma}(f_{2})(x)|+|I_{\gamma}[(b_{B(y,r)}-b)f_{2}](x)|.$$
Therefore,
\begin{align*}
\text{II}(y,r) &=|B(y,r)|^{\frac{1}{\beta}-\frac{1}{n}\sum_{i=1}^n\frac{1}{s_i}-\frac{1}{n}\sum_{i=1}^{n}\frac{1}{q_{i}}} \|[b,I_{\gamma}](f_{2})\chi_{B(y,r)}\|_{L^{\vec{q}}}\\
&\le|B(y,r)|^{\frac{1}{\beta}-\frac{1}{n}\sum_{i=1}^n\frac{1}{s_i}-\frac{1}{n}\sum_{i=1}^{n}\frac{1}{q_{i}}}\||b-b_{B(y,r)}| |\text{I}_{\gamma}(f_{2})|\chi_{B(y,r)}\|_{L^{\vec{q}}}\\
&+ \|I_{\gamma}[(b_{B(y,r) }-b)f_{2}]\chi_{B(y,r)}\|_{L^{\vec{q}}}\\
&\le|B(y,r)|^{\frac{1}{\beta}-\frac{1}{n}\sum_{i=1}^n\frac{1}{s_i}-\frac{1}{n}\sum_{i=1}^{n}\frac{1}{q_{i}}}
[\| |b-b_{B(y,r) }||I_{\gamma}(f_{2}) |\chi_{B(y,r)}\|_{L^{\vec{q}}}\\
&~~~+ |B(y,r)|^{\frac{1}{\beta}-\frac{1}{n}\sum_{i=1}^n\frac{1}{s_i}-\frac{1}{n}\sum_{i=1}^{n}\frac{1}{q_{i}}}
\|\text{I}_{\gamma}[(b_{B(y,r) }-b)f_{2}]\chi_{B(y,r)}\|_{L^{\vec{q}}}\\
&:=\text{II}_{1}(y,r)+\text{II}_{2}(y,r).
\end{align*}
According to (5.2) in the proof of theorem 2.9, we know
$$I_{\gamma}(f_{2})\lesssim \sum^{\infty}_{j=1} \displaystyle{\frac{1}{|2^{j+1}B(y,r)|^{1-\frac{\gamma}{n}}} }
\int _{2^{j+1}B(y,r)} |f(z)|dz.$$
By H\"oder's inequality and (ii) of Lemma 4.1,
\begin{align*}
\text{II}_{1}(y,r) &\lesssim|B(y,r)|^{\frac{1}{\beta}-\frac{1}{n}\sum_{i=1}^n\frac{1}{s_{i}}-\frac{1}{n}\sum_{i=1}^n\frac{1}{q_{i}}}
\sum^{\infty}_{j=1}\frac{\|(b-b_{B(y,r)})\chi_{B(y,r)}\|_{L^{\vec{q}}}}{|2^{j+1}| B(y,r)|^{1-\frac{\gamma}{n}}}\int _{2^{j+1} B(y,r)}|f(z)|dz\\
&\sim\sum^{\infty}_{j=1}|B(y,r)|^{\frac{1}{\beta}-\frac{1}{n}\sum_{i=1}^n\frac{1}{s_{i}}}
\displaystyle{\frac{1}{|2^{j+1} B(y,r)|^{1-\frac{\gamma}{n}}}}\int_{2^{j+1}B(y,r)}|f(z)|dz\cdot \displaystyle{\frac{\|(b-b_{B(y,r)})|\chi_{B(y,r)}\|_{L^{\vec{q}}}}
{\| \chi_{B(y,r)} \|_{L^{\vec{q}}}} }\\
&\sim \|b\|_{*} \sum^{\infty}_{j=1}
|B(y,r)|^{\frac{1}{\beta}-\frac{1}{n}\sum_{i=1}^n\frac{1}{s_{i}}}\frac{1}{|2^{j+1} B(y,r)|^{1-\frac{\gamma}{n}}}\int_{2^{j+1}B(y,r)}|f(z)|dz\\
&\leq \|b\|_{*} \sum^{\infty}_{j=1}
|B(y,r)|^{\frac{1}{\beta}-\frac{1}{n}\sum_{i=1}^n\frac{1}{s_{i}}}
|2^{j+1} B(y,r)|^{\frac{\gamma}{n}-1}
|2^{j+1} B(y,r)|^{\frac{1}{n}\sum_{i=1}^n\frac{1}{p_{i}}}
{\| f\chi_{2^{j+1}B(y,r)} \|_{L^{\vec{p}}}}\\
&=\|b\|_{*}\sum^{\infty}_{j=1}2^{-j(\frac{1}{\beta}-\frac{1}{n}\sum_{i=1}^n\frac{1}{s_{i}})} |2^{j+1}B(y,r)|^{\frac{1}{\alpha}-\frac{1}{n}\sum_{i=1}^n\frac{1}{s_{i}}-\frac{1}{n}\sum_{i=1}^n\frac{1}{p_{i}}}
{\|f\chi_{2^{j+1}B(y,r)} \|_{L^{\vec{p}}}}.
\end{align*}
Due to the assumption $\frac{1}{\beta}-\frac{1}{n}\sum_{i=1}^n\frac{1}{s_{i}}>0$,
$$\sum^{\infty}_{j=1}2^{-j(\frac{1}{\beta}-\frac{1}{n}\sum_{i=1}^n\frac{1}{s_{i}})}\sim 1.\eqno{(6.2)}$$
Thus,
$$\sup_{r>0} \| \text{I}(y,r)  \|_{L^{\vec{s}}(\R^n)}\lesssim \|b\|_{BMO} \|f\|_{(\L^{\vec{p}},L^{\vec{s}})^{\alpha}(\R^n)}.\eqno{(6.3)}$$

For the estimates of $\text{II}_{2}(y,r)$, we have
\begin{align*}
\text{II}_{2}(y,r)&=
|B(y,r)|^{\frac{1}{\beta}-\frac{1}{n}\sum_{i=1}^n\frac{1}{s_{i}}-\frac{1}{n}\sum_{i=1}^n\frac{1}{q_{i}}}
\| \text{I}_{\gamma} [(b_{B(y,r) }-b)f_{2}] \chi_{B(y,r)}\|_{L^{\vec{q}}}\\
&\le|B(y,r)|^{\frac{1}{\beta}-\frac{1}{n}\sum_{i=1}^n\frac{1}{s_{i}}-\frac{1}{n}\sum_{i=1}^n\frac{1}{q_{i}}}
\cdot\sum\limits_{j=1}\limits^{+\infty}\frac{1}{|2^{j+1}B(y,r)|^{1-\frac{\gamma}{n}}}\int_{2^{j+1}B(y,r)}|f(z)||b(z)-b_{B(y,r)}|dz
\cdot|B(y,r)|^{\frac{1}{n}\sum_{i=1}^n\frac{1}{q_{i}}}\\
&\le|B(y,r)|^{\frac{1}{\beta}-\frac{1}{n}\sum_{i=1}^n\frac{1}{s_{i}}}
\cdot\sum\limits_{j=1}\limits^{+\infty}\frac{1}{|2^{j+1}B(y,r)|^{1-\frac{\gamma}{n}}}\int_{2^{j+1}B(y,r)}|f(z)||b(z)-b_{2^{j+1}B(y,r)}|dz\\
&+|B(y,r)|^{\frac{1}{\beta}-\frac{1}{n}\sum_{i=1}^n\frac{1}{s_{i}}}
\cdot\sum\limits_{j=1}\limits^{+\infty}\frac{1}{|2^{j+1}B(y,r)|^{1-\frac{\gamma}{n}}}\int_{2^{j+1}B(y,r)}|f(z)|dz\cdot|b_{2^{j+1}B(y,r)}-b_{B(y,r)}|\\
&=:\text{II}_{21}(y,r)+\text{II}_{22}(y,r).
\end{align*}

To estimate $\text{II}_{21}(y,r)$, applying H\"older's inequality and the second part of Lemma 6.1, we can deduce that
\begin{align*}
\text{II}_{21}(y,r)&\leq|B(y,r)|^{\frac{1}{\beta}-\frac{1}{n}\sum_{i=1}^n\frac{1}{s_{i}}}
\cdot\sum\limits_{j=1}\limits^{+\infty}\frac{1}{|2^{j+1}B(y,r)|^{1-\frac{\gamma}{n}}}\|f\chi_{2^{j+1}B(y,r)}\|_{L^{\vec{p}}}
\|(b-b_{2^{j+1}B(y,r)})\chi_{2^{j+1}B(y,r)}\|_{L^{\vec{p}}}\\
&\sim\sum\limits_{j=1}\limits^{+\infty}|B(y,r)|^{\frac{1}{\beta}-\frac{1}{n}\sum_{i=1}^n\frac{1}{s_{i}}}
\cdot\frac{1}{|2^{j+1}B(y,r)|^{1-\frac{\gamma}{n}}}\cdot|2^{j+1}B(y,r)|^{\sum_{i=1}^{n}\frac{1}{p_{i}}}
\cdot\|f\chi_{2^{j+1}B(y,r)}\|_{L^{\vec{p}}}\\
&\times\|(b-b_{2^{j+1}B(y,r)})\chi_{2^{j+1}B(y,r)}\|_{L^{\vec{p}}}\cdot\|\chi_{2^{j+1}B(y,r)}\|_{L^{\vec{p}}}\\
&\sim\sum\limits_{j=1}\limits^{+\infty}2^{j(\frac{1}{\beta}-\frac{1}{n}\sum_{i=1}^n\frac{1}{s_{i}})}
\cdot|2^{j+1}B(y,r)|^{\frac{1}{\alpha}-\frac{1}{n}\sum_{i=1}^n\frac{1}{s_{i}}-\frac{1}{n}\sum_{i=1}^n\frac{1}{p_{i}}}
\cdot\|f\chi_{2^{j+1}B(y,r)}\|_{L^{\vec{p}}}\cdot\|b\|_{*}.
\end{align*}
By (4.2), we get
$$\sup_{r>0}\|\text{II}_{21}(y,r)\|_{L^{\vec{s}}}\lesssim\|b\|_{*}\cdot\|f\|_{(L^{\vec{p}},L^{\vec{s}})^{\alpha}}.\eqno{(6.4)}$$
Now, we estimate $\text{II}_{22}(y,r)$. An application of H\"older's inequality and first part of Lemma 4.1 gives us that
\begin{align*}
\text{II}_{22}(y,r)&\leq|B(y,r)|^{\frac{1}{\beta}-\frac{1}{n}\sum_{i=1}^n\frac{1}{s_{i}}}
\cdot\sum\limits_{j=1}\limits^{+\infty}\frac{1}{|2^{j+1}B(y,r)|^{1-\frac{\gamma}{n}}}\int_{2^{j+1}B(y,r)}|f(z)|dz \cdot|b_{2^{j+1}B(y,r)}-b_{B(y,r)}|\\
&\leq|B(y,r)|^{\frac{1}{\beta}-\frac{1}{n}\sum_{i=1}^n\frac{1}{s_{i}}}
\cdot\sum\limits_{j=1}\limits^{+\infty}\frac{1}{|2^{j+1}B(y,r)|^{1-\frac{\gamma}{n}}}\int_{2^{j+1}B(y,r)}|f(z)|dz\cdot(j+1)\|b\|_{BMO}\\
&\leq|B(y,r)|^{\frac{1}{\beta}-\frac{1}{n}\sum_{i=1}^n\frac{1}{s_{i}}}
\cdot\sum\limits_{j=1}\limits^{+\infty}\frac{j}{|2^{j+1}B(y,r)|^{1-\frac{\gamma}{n}}}\cdot\|f\chi_{2^{j+1}B(y,r)}\|_{L^{\vec{p}}}
\cdot|2^{j+1}B(y,r)|^{\frac{1}{n}\sum_{i=1}^n\frac{1}{p_{i}}}\cdot\|b\|_{BMO}\\
&\leq\sum\limits_{j=1}\limits^{+\infty}\frac{j}{2^{j(\frac{1}{\beta}-\frac{1}{n}\sum_{i=1}^n\frac{1}{s_{i}})}}
\cdot|2^{j+1}B(y,r)|^{\frac{1}{\alpha}-\frac{1}{n}\sum_{i=1}^n\frac{1}{s_{i}}-\frac{1}{n}\sum_{i=1}^n\frac{1}{p_{i}}}
\cdot\|f\chi_{2^{j+1}B(y,r)}\|_{L^{\vec{p}}}\cdot\|b\|_{BMO}.
\end{align*}
By (6.2), we get
$$\sup\limits_{r>0}\|\text{II}_{22}(y,r)\|_{L^{\vec{s}}}\le\|b\|_{BMO}\cdot\|f\|_{(L^{\vec{p}},L^{\vec{s}})^{\alpha}}.\eqno{(6.5)}$$
Combining (6.1), (6.3), (6.4), and (6.5), we conclude that
\begin{align*}
\|I_{\gamma}(f)\|_{(L^{\vec{p}},L^{\vec{s}})^{\alpha}}&\leq\sup\limits_{r>0}\|\text{I}(\cdot,r)\|_{L^{\vec{s}}}
+\sup\limits_{r>0}\|\text{II}(\cdot,r)\|_{L^{\vec{s}}}\\
&\leq\sup\limits_{r>0}\|\text{I}(\cdot,r)\|_{L^{\vec{s}}}+\sup\limits_{r>0}\|\text{II}_{1}(\cdot,r)\|_{L^{\vec{s}}}
+\sup\limits_{r>0}\|\text{II}_{2}(\cdot,r)\|_{L^{\vec{s}}}\\
&\leq\sup\limits_{r>0}\|\text{I}(\cdot,r)\|_{L^{\vec{s}}}+\sup\limits_{r>0}\|\text{II}_{1}(\cdot,r)\|_{L^{\vec{s}}}
+\sup\limits_{r>0}\|\text{II}_{21}(\cdot,r)\|_{L^{\vec{s}}}+\sup\limits_{r>0}\|\text{II}_{22}(\cdot,r)\|_{L^{\vec{s}}}\\
&\leq\|b\|_{*}\cdot\|f\|_{(L^{\vec{p}},L^{\vec{s}})^{\alpha}}.
\end{align*}
The proof is completed. $~~~~\blacksquare$

\section{Proof of Theorem 2.11}
\par
In this section, we prove the Theorem 2.11. Before that, we give an estimate of characteristic function on $(L^{\vec{p}},L^{\vec{s}})^\alpha(\R^n)$ and $\mathcal{H}(\vec{p}',\vec{s}\,',\alpha')(\R^n)$.

\textbf{Proposition 7.1.} Let $0\le\frac{1}{n}\sum_{i=1}^n\frac{1}{s_i}\leq\frac{1}{\alpha}\leq\frac{1}{n}\sum^{n}_{i=1}\frac{1}{p_{i}}<1$ and $\chi_{B(x_0,r_0)}$ is a characteristic function on $B(x_0,r_0)$. Then we have
$$\|\chi_{B(x_0,r_0)}\|_{(L^{\vec{p}},L^{\vec{s}})^\alpha}\lesssim r_0^{n/\alpha}\text{ and }\|\chi_{B(x_0,r_0)}\|_{\mathcal{H}(\vec{p}',\vec{s}\,',\alpha')}\lesssim r_0^{n/\alpha'}.$$

\textbf{Proof.} It is obviously that
$$\|\chi_{B(x_0,r_0)}\|_{(L^{\vec{p}},L^{\vec{s}})^\alpha} \sim\sup_{r>0}r^{\frac{n}{\alpha}-\sum^{n}_{i=1}\frac{1}{s_{i}}-\sum^{n}_{i=1}\frac{1}{p_{i}}}
\left\| \|\chi_{B(x_{0},r_0)}\chi_{B(x,r)}\|_{L^{\vec{p}}} \right\|_{L^{\vec{s}}}.$$
If $r>r_0$, then by $\frac{1}{\alpha}-\frac{1}{n}\sum^{n}_{i=1}\frac{1}{p_{i}}\le 0$
\begin{align*}
&\quad\sup_{r>r_0}r^{\frac{n}{\alpha}-\sum^{n}_{i=1}\frac{1}{s_{i}}-\sum^{n}_{i=1}\frac{1}{p_{i}}}
\left\| \|\chi_{B(x_{0},r_0)}\chi_{B(\cdot,r)}\|_{L^{\vec{p}}} \right\|_{L^{\vec{s}}}\\
&\le\sup_{r>r_0}r^{\frac{n}{\alpha}-\sum^{n}_{i=1}\frac{1}{s_{i}}-\sum^{n}_{i=1}\frac{1}{p_{i}}}
\left\| \|\chi_{B(x_{0},r_0)}\|_{L^{\vec{p}}}\cdot\chi_{B(x_0,r+r_0)}\right\|_{L^{\vec{s}}}\\
&\lesssim r_0^{\sum_{j=1}^n\frac{1}{p_j}}\sup_{r>r_0}r^{\frac{n}{\alpha}-\sum^{n}_{i=1}\frac{1}{s_{i}}-\sum^{n}_{i=1}\frac{1}{p_{i}}}
(r+r_0)^{\sum_{i=1}^n\frac{n}{s_i}}\\
&=r_0^{\sum_{j=1}^n\frac{1}{p_j}}\sup_{r>r_0}r^{\frac{n}{\alpha}-\sum^{n}_{i=1}\frac{1}{p_{i}}}
(1+\frac{r_0}{r})^{\sum^{n}_{i=1}\frac{1}{s_{i}}}\\
&\lesssim r_0^{\frac{n}{\alpha}}.
\end{align*}
For $r\le r_0$, by $\frac{1}{\alpha}-\frac{1}{n}\sum_{i=1}^n\frac{1}{s_i}\ge 0$ we have
\begin{align*}
&\quad\sup_{r\le r_0}r^{\frac{n}{\alpha}-\sum^{n}_{i=1}\frac{1}{s_{i}}-\sum^{n}_{i=1}\frac{1}{p_{i}}}
\left\| \|\chi_{B(x_{0},r_0)}\chi_{B(\cdot,r)}\|_{L^{\vec{p}}} \right\|_{L^{\vec{s}}}\\
&\le\sup_{r>r_0}r^{\frac{n}{\alpha}-\sum^{n}_{i=1}\frac{1}{s_{i}}-\sum^{n}_{i=1}\frac{1}{p_{i}}}
\left\| \|\chi_{B(\cdot,r)}\|_{L^{\vec{p}}}\cdot\chi_{B(x_0,r+r_0)}\right\|_{L^{\vec{s}}}\\
&\lesssim \sup_{r>0}r^{\frac{n}{\alpha}-\sum^{n}_{i=1}\frac{1}{s_{i}}-\sum^{n}_{i=1}\frac{1}{p_{i}}}r^{\sum_{j=1}^n\frac{1}{p_j}}
(r+r_0)^{\sum^{n}_{i=1}\frac{1}{s_{i}}}\\
&=\sup_{r>0}r^{\frac{n}{\alpha}-\sum^{n}_{i=1}\frac{1}{p_{i}}} (r+r_0)^{\sum^{n}_{i=1}\frac{1}{s_{i}}}\\
&\lesssim r_0^{\frac{n}{\alpha}}.
\end{align*}
Thus, $\|\chi_{B(x_0,r_0)}\|_{(L^{\vec{p}},L^{\vec{s}})^\alpha}\lesssim r_0^{n/\alpha}.$

Next, we show that $\|\chi_{B(x_0,r_0)}\|_{\mathcal{H}(\vec{p},\vec{s}\,',\alpha')}\lesssim r_0^{n/\alpha'}$. First, by the similar argument dilation operator of (4.3), let
$$\chi_{B(x_0,r_0)}=r^{\frac{n}{\alpha'}}\|\chi_{B(x_0/r,r_0/r)}\|_{\vec{p}',\vec{s}\,'}\cdot St_r^{\alpha'}(\|\chi_{B(x_0/r,r_0/r)}\|_{\vec{p}',\vec{s}\,'}^{-1}\chi_{B(x_0/r,r_0/r)}).$$
It is obvious that
$$\left\|\|\chi_{B(x_0/r,r_0/r)}\|_{\vec{p}',\vec{s}\,'}^{-1}\chi_{B(x_0/r,r_0/r)}\right\|_{\vec{p}',\vec{s}\,'}\le 1.$$
From Definition 2.6 and Proposition 2.5,
\begin{align*}
\|\chi_{B(x_0,r_0)}\|_{\mathcal{H}(\vec{p}\,',\vec{s}\,',\alpha')} &\le\sup_{r>0}r^{\frac{n}{\alpha'}}\|\chi_{B(x_0/r,r_0/r)}\|_{\vec{p}',\vec{s}\,'}\\
&\lesssim \sup_{r>0}r^{\frac{n}{\alpha'}}\left\|\|\chi_{B(x_0/r,r_0/r)}\chi_{B(\cdot,1)}\|_{L^{\vec{p}'}} \right\|_{L^{\vec{s}\,'}}.
\end{align*}
Using the same argument of the proof of $\|\chi_{B(x_0,r_0)}\|_{(L^{\vec{p}},L^{\vec{s}})^\alpha}\lesssim r_0^{n/\alpha}$ with $r_0/r>1$ and $r_0/r\le 1$, we have
$$\|\chi_{B(x_0,r_0)}\|_{\mathcal{H}(\vec{p}\,',\vec{s}\,',\alpha')} \lesssim\sup_{r>0}r^{\frac{n}{\alpha'}}\left\|\|\chi_{B(x_0/r,r_0/r)}\chi_{B(\cdot,1)}\|_{L^{\vec{p}'}} \right\|_{L^{\vec{s}\,'}}\lesssim r_0^{n/\alpha'}.$$
The proof is completed. $~~~~\blacksquare$

Now, let us prove Theorem 2.11.

\textbf{Proof of Theorem 2.11.}  Assume that $[b,I_\alpha]$ is bounded from $(L^{\vec{p}},L^{\vec{s}})^{\alpha}(\R^n)$ to $(L^{\vec{q}},L^{\vec{s}})^{\beta}(\R^n)$. We use the same method as Janson \cite{32}. Choose $0\neq z_0\in\mathbb{R}^n$ such that $0\notin B(z_0,2)$. Then for $x\in B(z_0,2)$, $|x|^{n-\alpha}\in C^{\infty}(B(z_0,2))$. Hence, $|x|^{n-\alpha}$ can be written as the absolutely convergent Fourier series:
$$|x|^{n-\alpha}\chi_{B(z_0,2)}(x)=\sum_{m\in \mathbb{Z}^n}a_me^{2im\cdot x}\chi_{B(z_0,2)}(x)$$
with $\sum_{m\in \mathbb{Z}^n}|a_m|<\infty$.

For any $x_0\in\mathbb{R}^n$ and $t>0$, let $B=B(x_0,t)$ and $B_{z_0}=B(x_0+z_0t,t)$. Let $s(x)=\overline{sgn(\int_{B_{z_0}}(b(x)-b(y))dy)}$. Then
\begin{align*}
\frac{1}{|B|}\int_B|b(x)-b_{B_{z_0}}|
=\frac{1}{|B|}\frac{1}{|B_{z_0}|}\int_B\int_{B_{z_0}}s(x)(b(x)-b(y))dydx.
\end{align*}
If $x\in B$ and $y\in B_{z_0}$, then $\frac{y-x}{t}\in B(z_0,2)$. Thereby,
\begin{align*}
&\quad\frac{1}{|B|}\int_B|b(x)-b_{B_{z_0}}|\\
&=t^{-n-\gamma}\int_B\int_{B_{z_0}}s(x)(b(x)-b(y))|x-y|^{\alpha-n}\left(\frac{|x-y|}{t}\right)^{n-\gamma}dydx\\
&=t^{-n-\gamma}\sum_{m\in\mathbb{Z}^n}a_m\int_B\int_{B_{z_0}}s(x)(b(x)-b(y))|x-y|^{\gamma-n}e^{-2im\cdot \frac{y}{t}}dy\times e^{2im\cdot \frac{x}{t}}dx\\
&=t^{-n-\gamma}\sum_{m\in\mathbb{Z}^n}a_m\int_B[b,I_\gamma](e^{-2im\cdot \frac{\cdot}{t}}\chi_{B_{z_0}})(x)\times s(x)e^{2im\cdot \frac{x}{t}}dx.
\end{align*}
By (2.4) and Proposition 2.5,
$$\frac{1}{|B|}\int_B|b(x)-b_{B_{z_0}}|
\lesssim t^{-n-\gamma}\sum_{m\in\mathbb{Z}^n}a_m\left\|[b,I_\gamma](e^{-2im\cdot \frac{\cdot}{t}}\chi_{B_{z_0}})\right\|_{(L^{\vec{q}},L^{\vec{s}})^{\beta}} \left\|s\cdot e^{-2im\cdot \frac{\cdot}{t}}\chi_B\right\|_{\mathcal{H}(\vec{q}',\vec{s}\,',\beta')}.$$
It is easy to calculate
$$\left\|s\cdot e^{-2im\cdot \frac{\cdot}{t}}\chi_B\right\|_{\mathcal{H}(\vec{q}',\vec{s}\,',\beta')}\lesssim t^{n/\beta'}.$$
Hence,
$$\frac{1}{|B|}\int_B|b(x)-b_{B_{z_0}}|\lesssim t^{-n-\gamma+n/\beta'}\sum_{m\in\mathbb{Z}^n}a_m\left\|[b,I_\gamma](e^{-2im\cdot \frac{\cdot}{t}}\chi_{B_{z_0}})\right\|_{(L^{\vec{q}},L^{\vec{s}})^{\beta}}.$$
According to the hypothesis
\begin{align*}
\frac{1}{|B|}\int_B|b(x)-b_{B_{z_0}}|&\lesssim t^{-n-\gamma+n/\beta'}\left\|[b,I_\gamma]\right\|\sum_{m\in\mathbb{Z}^n}a_m\left\|e^{-2im\cdot \frac{\cdot}{t}}\chi_{B_{z_0}}\right\|_{(L^{\vec{p}},L^{\vec{s}})^\alpha}\\
&\le t^{-n-\gamma+n/\beta'+n/\alpha}\left\|[b,I_\gamma]\right\|\sum_{m\in\mathbb{Z}^n}a_m\\
&\lesssim \left\|[b,I_\gamma]\right\|.
\end{align*}
Thus, we have
$$\frac{1}{|B|}\int_B|b(x)-b(y)|dx\le\frac{2}{|Q|}\int_B|b(x)-b_{B_{z_0}}|dx\lesssim \left\|[b,I_\gamma]\right\|.$$
This prove $b\in BMO(\mathbb{R}^n)$. $~~~~\blacksquare$\\

\hspace*{-0.6cm}\textbf{\bf Acknowledgments}\\
The authors would like to express their thanks to the referees for valuable advice regarding previous version of this paper. This project is supported by the National Natural Science Foundation of China (Grant No.12061069).\\

\end{document}